\documentclass[a4paper]{amsart}\usepackage{amsfonts,amssymb,bbm,graphicx} 
\usepackage{enumerate}  

\usepackage{amsthm}

\def\A{\Bbb{A}}\def\C{\Bbb{C}}\def\k{\mathbbm{k}}\def\K{\mathbb{K}}
\def\bk{{\bar{\k}}}\def\N{\Bbb{N}}
\def\R{\Bbb{R}}\def\Z{\Bbb{Z}}

\def\di{\partial}

\def\suml{\sum\limits}
\def\capl{\mathop\cap\limits}
\def\cupl{\mathop\cup\limits}

\def\wedgel{\mathop\wedge\limits}

\newcommand{\quot}[2]{{\left.\raisebox{1.6ex}{$#1$}\!\!\!\!\!{\scalebox{2}{\ensuremath\diagup}}
\!\!\!\!\!\raisebox{-1ex}{$#2$}\right.}}
\newcommand{\quots}[2]{{\footnotesize\left.\raisebox{0.4ex}{$#1$}\! / \!\raisebox{-0.4ex}{$#2$}\right.}}

\def\tA{\tilde{A}}

\def\hR{{\widehat{R}}}

\def\al{\alpha}\def\be{\beta}

\def\la{\lambda}\def\Om{{\Omega}}\def\Si{\Sigma}

\def\cA{\mathcal A}

\def\cD{\mathcal D}
\def\cG{\mathcal G}\def\cK{{\mathcal K}}\def\cR{\mathfrak{R}}
\def\cm{{\frak m}}\def\cp{{\frak p}}\def\cq{{\frak q}}

\def\uf{\underline{f}}

\def\ux{\underline{x}}

\def\one{{1\hspace{-0.1cm}\rm I}}\def\zero{\mathbb{O}}

\newcommand{\ber}{\begin{array}{l}}\newcommand{\eer}{\end{array}}
\newcommand{\bpm}{\begin{pmatrix}}\newcommand{\epm}{\end{pmatrix}}
\newcommand{\bbm}{\begin{bmatrix}}\newcommand{\ebm}{\end{bmatrix}}
\newcommand{\bM}{\begin{matrix}}\newcommand{\eM}{\end{matrix}}
\newcommand{\bee}{\begin{enumerate}}\newcommand{\eee}{\end{enumerate}}
\newcommand{\bei}{\begin{itemize}}\newcommand{\eei}{\end{itemize}}

\def\Span{{\rm Span}}

\def\sset{\subset}\def\sseteq{\subseteq}\def\ssetneq{\subsetneq}\def\smin{\setminus}

\def\Mat{Mat_{m\times n}(R)}\def\Matm{Mat_{m\times n}(\cm)}

\newtheorem{Lemma}{Lemma}[section]\newcommand{\bel}{\begin{Lemma}}\newcommand{\eel}{\end{Lemma}}

\newtheorem{Theorem}[Lemma]{Theorem}\newcommand{\bthe}{\begin{Theorem}}\newcommand{\ethe}{\end{Theorem}}
\newtheorem{Proposition}[Lemma]{Proposition}\newcommand{\bprop}{\begin{Proposition}}\newcommand{\eprop}{\end{Proposition}}
\newtheorem{Corollary}[Lemma]{Corollary}\newcommand{\bcor}{\begin{Corollary}}\newcommand{\ecor}{\end{Corollary}}
\newtheorem{Definition}[Lemma]{Definition}\newcommand{\bed}{\begin{Definition}}\newcommand{\eed}{\end{Definition}}
\newtheorem{Definition-Proposition}[Lemma]{Definition-Proposition}

\def\bpr{{\em Proof.\ }}
\newcommand{\epr}{{\hfill\ensuremath\blacksquare}\\}
\newtheorem{Remark}[Lemma]{Remark}\newcommand{\beR}{\begin{Remark}\rm}\newcommand{\eeR}{\end{Remark}}
\newtheorem{Example}[Lemma]{Example}\newcommand{\bex}{\begin{Example}\rm}\newcommand{\eex}{\end{Example}}
\newtheorem{Problem}[Lemma]{Problem}\newcommand{\bprob}{\begin{Problem}\rm}\newcommand{\eprob}{\end{Problem}}
\newtheorem{Properties}[Lemma]{Properties}

\newcommand{\bet}{\begin{tabular}{cccccccc}}\newcommand{\eet}{\end{tabular}}
\newcommand{\beq}{\begin{equation}}\newcommand{\eeq}{\end{equation}}

\newcommand{\bin}[2]{\binom{#1}{#2}}
\newcommand\isom{\xrightarrow{\,\smash{\raisebox{-0.65ex}{\ensuremath{\scriptstyle\sim}}}\,}}

\title[]{G\MakeLowercase{roup actions on matrices over local rings.}
\\A\MakeLowercase{nnihilators of} $T^1$-\MakeLowercase{modules for the groups} $\cG_{lr}$, $\cG_{congr}$.
}
\author{D\MakeLowercase{mitry} K\MakeLowercase{erner}} 
\address{Department of Mathematics, Ben Gurion University of the Negev, P.O.B. 653, Be'er Sheva 84105, Israel.}
\email{dmitry.kerner@gmail.com}
\date{\today}
\thanks{D.K. was supported by the grant No. 844/14  of Israel Science Foundation. 
}

\begin{document}\maketitle
\begin{abstract}
We consider matrices with entries in a local ring, $A\in \Mat$. Fix a group action, $G\circlearrowright\Mat$, and a subset
 of allowed deformations, $\Si\sseteq\Mat$. The traditional objects of study in Singularity Theory and Algebraic Geometry are
 the tangent spaces $T_{(\Si,A)}$, $T_{(GA,A)}$, and their quotient, the tangent module to
  the miniversal deformation, $T^1_{(\Si,G,A)}=\quots{T_{(\Si,A)}}{T_{(GA,A)}}$.

    This module plays the key role in various deformation problems, e.g., deformations of maps, of modules, of (skew-)symmetric forms.
 In particular, the first question is to determine the support/annihilator of this tangent module.
 In \cite{Belitski-Kerner1} we have studied this tangent module for various $R$-linear group actions $G\circlearrowright\Mat$

In the current work we study the support of the module  $T^1_{(\Si,G,A)}$ for group actions that involve automorphisms of the ring.
 (Geometrically, these are group actions that involve the  local coordinate changes.)

We obtain various bounds on localizations of $T^1_{(\Si,G,A)}$ and compute the radical of the annihilator of $T^1_{(\Si,G,A)}$,
 i.e., the set-theoretic support. This brings the definition of an (apparently new) type  of singular locus,
 the   ``essential singular locus" of a map/sub-scheme.
 It reflects the ``unexpected" singularities of a subscheme, ignoring those imposed by the singularities of the ambient space.
 Unlike the classical singular locus (defined by a Fitting ideal of the module of differentials) the essential is defined by the annihilator ideal of
  the module of derivations.
\end{abstract}\setcounter{secnumdepth}{6} \setcounter{tocdepth}{1}
\tableofcontents

\section{Introduction}
\subsection{Setup}
Let $(R,\cm)$ be a commutative Noetherian local ring over a field $\k$ of zero characteristic.
 (The typical cases are $R=\quots{\k[[x]]}{J}$ or $R=\quots{\k\{x\}}{J}$ where $\k$ is a complete normed field. Here $x=(x_1,\dots,x_p)$.)
 Denote by $\Mat$ the $R$-module of $m\times n$ matrices. We always assume $m\le n$.
Various groups act on $\Mat$ and show up in various areas.
\bex\label{Ex.Intro.Typical.Groups}
\bee[i.]
\item
 If the action $G\circlearrowright\Mat$ is $R$-linear and preserves the subset of degenerate matrices then $G$ is contained in the group of left-right multiplications,
   $G_{lr}:=GL(m,R)\times GL(n,R)$.
  (See \cite[\S3.6]{Belitski-Kerner1} for the precise statement.)
 Matrices considered up to $GL(n,R)$-transformations correspond to  embedded modules,  $Im(A)\sset R^{\oplus m}$. Matrices considered
   up to  $GL(m,R)\times GL(n,R)$-transformations correspond to  non-embedded modules,  $Coker(A)=\quots{R^{\oplus m}}{Im(A)}$.
\item The group of $\k$-linear ring automorphisms, $Aut_\k(R)$, acts on matrices entry-wise.
 Geometrically they are the local coordinate changes on $Spec(R)$. In Singularity Theory this group is known as the right equivalence, $\cR$.

\item Accordingly one considers the semi-direct products, $\cG_{r}:=GL(n,R)\rtimes Aut_\k(R)$, $\cG_{lr}:=GL(m,R)\times GL(n,R)\rtimes Aut_\k(R)$.
  The action on the modules $Im(A)$, $Coker(A)$ is by the base change, $Coker(A)\to \phi(R)\underset{R}{\otimes}Coker(A)$.
 For one-row matrices, $m=1$, the orbits of $\cG_{lr}$, $\cG_r$ coincide. In Singularity Theory this group is known as the contact equivalence of maps, $\cK$.

\item The congruence, $G_{congr}=GL(m,R)\circlearrowright Mat_{m\times m}(R)$, acts by $A\to UAU^t$.
 Matrices considered up to the congruence correspond to   bilinear/symmetric/skew-symmetric forms.
 Accordingly one considers $\cG_{congr}:=G_{congr}\rtimes Aut_\k(R)$.
\eee
\eex

The traditional approach of deformation theory is to study the tangent space to the miniversal deformation.
 In our case this is the Tjurina module, $T^1_{(\Si,G,A)}:=\quots{T_{(\Si,A)}}{T_{(GA,A)}}$.
Here $\Si$ is one of $\Mat$, $Mat^{sym}_{m\times m}(R)$,  $Mat^{skew-sym}_{m\times m}(R)$, while $T_{(GA,A)}:=T_{(G,\one)}A$ is the image
 tangent space to the orbit.
 (For our particular cases these are defined in \S \ref{Sec.Background.Tangent.Spaces}.)

 In the simplest case of ``functions", $m=1=n$, $\Si=R$, we get the Milnor algebra, $T^1_{(R,Aut_\k(R),A)}$, and the Tjurina algebra,
 $T^1_{(R,\cG_{lr},A)}$.

\subsection{}

The module $T^1_{(\Si,G,A)}$ is complicated. It is usually a torsion over $R$, while over $\quots{R}{Ann(T^1_{(\Si,G,A)})}$ it is far from being
 free and usually of high rank.
 The complexity and importance can be appreciated by the following particular cases.
\bex
\bee[\bf i.]
\item
For $\Si=\Mat$, $G=GL(m,R)\times GL(n,R)$ and $A\in Mat_{m\times n}(\cm)$ there holds:
\[T^1_{(\Si,G,A)}=Ext^1_R\Big(Coker(A),Coker(A)\Big).\]
\item
Let $m=1<n$, and identify $Mat_{1\times n}(R)\approx R^n$, then $T^1_{(R^n,\cG_{lr},A)}$ is the classical Tjurina module
  of the map  $Spec(R)\stackrel{A}{\to}(\k^n,0)$.
  This module is among the cornerstones of the Singularity Theory, see chapter 4 of  \cite{Looijenga}.
   Its structure is rich and is not completely understood. For $R$ regular $T^1$ defines the singular locus of the map $A$.
\eee
In both cases one cannot present  $T^1_{(\Si,G,A)}$ (or its annihilator) in any simple/more explicit form.
\eex

\

In this paper we study the annihilator/support of $T^1_{(\Si,G,A)}$ for the actions $\cG_{lr}$, $\cG_{congr}$ of example \ref{Ex.Intro.Typical.Groups}.
 This information is
   needed e.g., for the studies of determinantal singularities, see   \cite{Bruce-Goryun-Zakal02}, \cite{Bruce-Tari04}, \cite{Bruce2003},
    \cite{Fr�hbis-Kr�ger.99},
     \cite{Fr�hbis-Kr�ger.15}, \cite{Fr�hbis-Kr�ger.18}, \cite{Damon-Pike},  \cite{Ahmed-Ruas}, \cite{N.B.-O.O.-T.13}, \cite{N.B.-O.O.-T.18}.
      In particular, if the annihilator $Ann(T^1_{(\Si,G,A)})$ contains a power of
  the maximal ideal (and the ring $R$ is henselian), then the studied object is ``finitely determined", \cite{Belitski-Kerner}, \cite{BGK}.
 Equivalently, the object is ``algebraizable in families"  \cite[\S3.12]{Belitski-Kerner1}).
 This is the first step in establishing  the finite dimensionality of the miniversal deformation,
  and then possible classification of simple/unimodal/\dots    singularity types.

The related questions of finite determinacy were studied in some particular cases, over the rings $\k[[\ux]]$, $\C\{\ux\}$, \cite{Greuel-Pham.17.a},
 \cite{Greuel-Pham.17.b}, \cite{Greuel-Pham.18}, but the approach
 was mostly algorithmic, translating the question into the case-by-case tasks for computer packages.
 Unlike the previous studies, we work in the generality of local Noetherian rings (without any regularity assumptions) and give explicit
  criteria, applicable without computer help.

\

     The bounds on $Ann(T^1_{(\Si,G,A)})$ in theorems \ref{Thm.Results.T1.Grl.coord.changes},
\ref{Thm.Results.T1.for.Congruence}  are somewhat involved.
 This is not a surprise, noticing the ``complexity of the problem". Recall that $T^1_{(\Si,G,A)}$ encodes many of the singularity properties of the object,
  and the bounds clearly show the singularity invariants of the matrices.
 On the other hand, the bounds admit transparent geometric interpretation, in terms of certain degeneracy loci. We emphasize also that the
 formulae are ``computationally simple", and admit direct computer implementation.

 In \cite{Belitski-Kerner2}  we apply these bounds  to establish strong results on finite determinacy/algebraizability/properties of the miniversal deformation
  of matrices  under various group actions. This gives criteria of algebraizability of  embedded modules, (skew-)symmetric forms,
   complexes of modules, etc. (The statements in \cite{Belitski-Kerner2} are of the following type.
   Let $M_\hR=coker(A)$ be a finitely generated module over a complete Noetherian local ring $\hR$, over $\k$. If $\sqrt{Ann(T^1_{(\Si,\cG_{lr},A)})}=\cm$ then
   there exists a finitely generated $\k$ subalgebra $R\sset \hR$, whose completion is $\hR$, and a finitely generated module $M_R$ such that $M_\hR=\hR\otimes M_R$.)

\

Some of our statements hold in the generality of commutative unital rings, for some others it is enough to assume that $R$ is local
 Noetherian (but not necessarily over a field of zero characteristic). Yet, we restrict to ``$\k$ is a field of zero characteristic",
  to avoid   various pathologies of the modules of deivations $Der_\k(R)$, and differentials, $\Om^1_{R/\k}$.
   In fact, our work has originated from the classical Singularity
   Theory, over  the classical rings  like $\quots{\k[[\ux]]}{J}$, $\quots{\C\{\ux\}}{J}$.  Already for these rings our results are new.

\subsection{The content and the structure of the paper}
\bee[i.]
\item The bounds on $Ann(T^1_{(\Si,G,A)})$ are stated in \S\ref{Sec.Results}.
 They are expressed in terms of the ideal $Sing_r(J)\sset R$, the ``essential singular locus" of the subscheme $V(J)\sset Spec(R)$, see below.
 This essential singular locus measures  the singularities of the determinantal strata, $Sing_r(I_j(A))$, here $r$ is the expected grade of $I_j(A)$.
 Recall that the scheme $V(I_j(A))\sset Spec(R)$ is always singular along $V(I_{j-1}(A))$.
  The module $T^1_{(\Si,G,A)}$ is supported exactly at the ``unexpected" singular points, $V(Sing_r(I_j(A)))\smin V(I_{j-1}(A))$.

 While the bounds are algebraic, we give transparent geometric interpretations. These bounds are extensively used in \cite{Belitski-Kerner2}.

 \item
In \S\ref{Sec.Preparations} we prepare the tools.
  Trying to be understandable by non-experts in Commutative Algebra we emphasize the geometric
 meaning of various notions and recall some standard results.
\bee
\item In \S\ref{Sec.Background.Tangent.Spaces}   we describe the tangent spaces to the
 group orbits, $T_{(G,\one)}A$.
\item  In \S\ref{Sec.Background.Localization} we establish the needed properties of localization.
\item  Sections \ref{Sec.Background.Determinantal.and.Pfaffians}, \ref{Sec.Background.Ann.Coker}, \ref{Sec.Background.Generaliz.Ann.Coker} are about the properties of the
 determinantal ideals, $\{I_j(A)\}$, Pfaffian ideals, $\{Pf_j(A)\}$, the  annihilator-of-cokernel, $Ann.Coker(A)$,
and its generalizations, $\{Ann.Coker_j(A)\}$.
\item In \S\ref{Sec.Background.Sing(J)} we study the essential singular locus, $Sing_r(J)\sset R$.
 (Here $r$ is the expected grade of $J\sset R$.)

Recall that the classical singular locus is  defined by the Fitting ideal of the module of differentials, $Fitt_{dim(\quots{R}{J})} \Om^1_{R/J}\sset R$.
 The essential singular locus is defined using the module of derivations, $Der_\k(R)$, and the annihilator scheme structure.
 Unlike the classical singularity locus, $Sing_r(J)$ measures only the ``unexpected" singularities of $V(J)\sset Spec(R)$ and
 often ignores the singularities of the ambient space, $Spec(R)$.

 If $R$ is a regular local ring and the ideal $J$ is pure, of grade $r$, we get the classical singular locus
  (but with annihilator rather than Fitting ideal scheme structure):
  $  Sing_r(J)=Ann_r\Om^1_{R/J}$. But for non-regular rings or non-pure ideals, one
   has usually:
\beq
     Sing_r(J)\supsetneq Ann_r\Om^1_{R/J}\supsetneq Fitt_{dim(\quots{R}{J})} \Om^1_{R/J}.
\eeq
\item In \S\ref{Sec.Background.Invariance.Annihilator} we prove that the ideal $Ann(T^1_{(\Si,G,A)})$ is invariant under the action of some elements of $\cG_{lr}$.
 This fact is used repeatedly in \S\ref{Sec.Proofs}.

\eee
\item
In \S\ref{Sec.Proofs} we prove the statements of \S2.
 The proofs go by checking the support of  $T^1_{(\Si,G,A)}$ ``pointwise", i.e., by localizations at prime ideals.

\eee

\subsection{Notations and conventions}\label{Sec.Results.Notations} 

\bee[\bf i.]
\item
The ideal quotient is $I:J=\{f\in R|\ f\cdot J\sseteq I\}$.
\item
The saturation of $I\sset R$ by $J\sset R$ is the ideal $Sat_J(I):=\suml_{k=1}^\infty I:J^k$.
 (We prefer not to use the standard notation $I:J^\infty$ to avoid any confusion with the ideal $J^\infty:=\cap_k J^k$.)
 Geometrically one
erases the subscheme $V(\sqrt{J})\sset Spec(R)$ and then takes the
Zariski closure, $V(Sat_J(I))=\overline{V(I)\smin V(\sqrt{J})}$.

We often use the relation $\sqrt{Sat_J(I)}=\sqrt{I}:J$, see lemma \ref{Thm.Background.Saturation}.
\item Localization at a prime ideal, $R\to R_\cp$, induces the natural map $\Mat\to Mat_{m\times n}(R_\cp)$. We denote the image of $A\in \Mat$
 by $A_\cp\in Mat_{m\times n}(R_\cp)$.

\item Suppose an $R$-module $M$ is minimaly generated by $m$ elements.
 The annihilator ideal  $Ann(M)$  is a refined version of the Fitting ideal, $Fitt_0(M)$. Similarly, the $j$'th annihilator, $Ann_j(M)$,
 is the refinement of the ideal $Fitt_{m-j}(M)$. Choosing a particular presentation matrix of a module, $M=Coker(A)$,
  we get the determinantal ideals, $I_j(A)=Fitt_{m-j}(Coker(A))$, and their refined versions, $\{Ann.Coker_j(A):=Ann_{j}(Coker(A))\}$.
 See \S\ref{Sec.Background.Ann.Coker}, \S\ref{Sec.Background.Generaliz.Ann.Coker} for the  definitions and properties.

\item
Let $Der_\k(R)$ be the module of ($\k$-linear) derivations of $R$.
The derivations act on matrices entrywise, for any $\cD\in Der_\k(R)$ one has $\cD(A)\in\Mat$. By applying the whole module $Der_\k(R)$
 we get the submodule $Der_\k(R)(A)\sseteq \Mat$. Similarly, for an ideal $J\sset R$ one gets the ideal $Der_\k(R)(J)\sset R$.

Sometimes we need only the submodule $Der_\k(R,\cm)$, the derivations sending $R$ to $\cm$. Accordingly we have
 $Der_\k(R,\cm)(J)\sset R$  and  $Der_\k(R,\cm)(A)\sseteq Mat_{m\times n}(\cm)$.

\item
Fix an ideal, $J\sseteq R$, and a number $r\in \N$, which is often the expected height/grade of $r$.
 Assume $J$ is finitely generated, choose any system of
generators, $J=\langle f_1,\dots,f_N\rangle$, write them as a column, $\uf$. Applying
the derivations of $R$ to this column we get the submodule $Der_\k(R)(\uf)\sseteq R^{\oplus N}$.

The essential singular locus of  $J$ is defined as
\beq\label{Eq.Intro.def.of.Sing(J)}
Sing_r(J):=\Bigg\{\ber
Ann_r\quot{R^{\oplus N}}{J\cdot R^{\oplus
N}+Der_\k(R)(\uf)},\quad \text{for }r\le N;
\\J, \quad \text{for }r> N.
\eer
\eeq
(Here $Ann_r$ is the $r$'th annihilator ideal, a refinement of the $r$'th determinantal ideal, see \S\ref{Sec.Background.Generaliz.Ann.Coker}.)

\

The typical context for $Sing_r(J)$ is the determinantal ideals, $J=I_{j+1}(A)$. Then $r$ is taken
 as the expected height/grade:
\bee
\item for $A\in \Matm$ one takes $r=(m-j)(n-j)$;
\item for  $A\in Mat^{sym}_{m\times m}(\cm)$ one takes $r=\bin{m+1-j}{2}$
\item for $A\in Mat^{skew-sym}_{m\times m}(\cm)$ and $j$-even one takes $r=\bin{m-j}{2}$;
 \item for  $A\in Mat^{skew-sym}_{m\times m}(\cm)$ and $j$-odd one takes $r=\bin{m-j+1}{2}$
\eee
Note that we do not take here the usual  $min\big((m-j)(n-j),dim(R)\big)$ or  $min\big((m-j)(n-j),depth(R)\big)$.

\

Sometimes we use the $\cm$-singular locus, $Sing^{(\cm)}_r(J)$, with $Der_\k(R,\cm)$ instead of $Der_\k(R)$.
\eee

\subsection{Acknowledgement}
Many thanks to A.Fernandez-Boix and G. Belitskii for useful advices.

\section{The main results}\label{Sec.Results}

\subsection{}
The action $Aut_\k(R)\circlearrowright \Si:=\Mat$ does not involve any matrix structure and the presentation of
  the tangent space $T_{(Aut_k(R)A,A)}$ (\S\ref{Sec.Background.Tangent.Spaces})
gives the obvious
\beq
Ann\big(T^1_{(\Si,Aut_\k(R),A)}\big)=Ann\big(\quots{\Si}{Der_\k(R,\cm)(A)}\big).
\eeq
For the regular rings like $\k[[\ux]]$, $\C\{\ux\}$ this recovers the classically studied cases in Singularity Theory, e.g.
 the ``local ring version" of the Milnor algebra for $m=n=1$. For the non-regular rings $\quots{\k[[\ux]]}{J}$, $\quots{\C\{\ux\}}{J}$
  see e.g., \cite{Bruce-Roberts}.

We remark that for $(m,n)\neq(1,1)$ and $A\in Mat_{m\times n}(\cm^2)$ the ideal $Ann\big(\quots{\Si}{Der_\k(R,\cm)(A)}\big)$ does not contain any power of $\cm$.
 (See e.g., proposition 5.7 of \cite{BGK}.)

\subsection{The action $\cG_{lr}\circlearrowright\Mat$ }\label{Sec.Results.Criteria.G.coord.with.changes}
\bthe\label{Thm.Results.T1.Grl.coord.changes}
Let $R$ be a Noetherian local ring over a field $\k$ of zero characteristic,  $A\in \Si:=\Mat$.
\bee[1.]
\item If $m=1$ then $Ann\big(T^1_{(\Si,\cG_{lr},A)}\big)=Sing^{(\cm)}_{n}(I_1(A))\supseteq I_1(A)+Ann\big(\quots{\Si}{Der_\k(R,\cm)(A)}\big)$.
\item For any prime ideal  satisfying $\cm\supsetneq\cp\supseteq I_m(A)$ but $\cp\not\supseteq I_{m-1}(A)$, the
 localizations of ideals at $\cp$ satisfy:
\[Ann(T^1_{(\Si,\cG_{lr},A)})_{\cp}=Sing_{n-m+1}(I_m(A))_{\cp}.\]
\item 
$Ann.Coker(A)+Ann\big(\quots{\Si}{Der_\k(R,\cm)(A)}\big)\sseteq Ann(T^1_{(\Si,\cG_{lr},A)})\sseteq
\capl^{m-1}_{j=0}Sat_{I_j(A)}\Big(Sing^{(\cm)}_{(m-j)(n-j)}\big(I_{j+1}(A)\big)\Big)$.
\item Suppose  for some $j$ the prime decomposition is $\sqrt{I_{j+1}(A)}=(\cap_\al \cp_\al)\cap(\cap_\be \cq_\be)$,
 where $grade(\cp_\al)<(m-j)(n-j)$ and $grade(\cq_\be)=(m-j)(n-j)$.
 Take the corresponding localizations, $R\stackrel{i_\al}{\to}R_{\cp_\al}$, and suppose the rings  $\{R_{\cp_\al}\}$ are regular.   Then
$Ann(T^1_{(\cG_{lr},\Si,A)})\sseteq Sat_{I_j(A)}\Big(\cap_\al i^{-1}_\al(I_{j+1}(A)_{\cp_\al})\Big)$.
\item $\sqrt{Ann(T^1_{(\Si,\cG_{lr},A)})}=
\capl^{m-1}_{j=0}\Big(\sqrt{Sing^{(\cm)}_{(m-j)(n-j)}\big(I_{j+1}(A)\big)}:I_{j}(A)\Big)$.
\eee
\ethe
Below are some remarks and  geometric interpretations (when the base field is algebraically closed, $\k=\bk$).
\bee[\bf {Part} 1.]
\item Here the matrix structure plays no role and  the $\cG_{lr}$-action induces the classical contact equivalence of maps, $\cK$.
 For $m=1=n$ (re)denote $A$ by $f$, then $Ann\big(T^1_{(\Si,\cG_{lr},f)}\big)=(f)+Der_\k(R,\cm)(f)$. This is the ``local ring version"
  of the classical Tjurina ideal of a function.

  More generally, for  $m=1\le n$, the ideal $Sing_n(I_1(A))$ defines the essential singular locus of the subscheme $A^{-1}(0)\sset Spec(R)$
   or of the map $Spec(R)\stackrel{A}{\to}(\k^n,0)$. Recall that we use  derivations instead of  differentials
    and  the annihilator scheme structure, rather than the Fitting scheme structure.
\item $T^1_{(\Si,\cG_{lr},A)}$ is supported on the essential singular locus of $V(I_m(A))$, with the locus $V(I_{m-1}(A))$ erased:
\[
Supp(T^1_{(\Si,\cG_{lr},A)})\smin\Big(Supp(T^1_{(\Si,\cG_{lr},A)})\cap V(I_{m-1}(A))\Big)=Sing_{n-m+1}\Big(V(I_m(A))\Big)\smin V(I_{m-1}(A))\sset Spec(R).
\]
Recall that $A$ is ``infinitesimally $G_{lr}$-stable", i.e. $T_{(G_{lr},A)}=\Mat$, at the points of $Spec(R)\smin V(I_m(A))$, where $A$ is of the full rank.
 Part 2 gives: $A$ is ``infinitesimally $\cG_{lr}$-stable" at the points of $Spec(R)\smin V(Sing_{n-m+1}(I_m(A)))$.

\item The embedding $Ann(T^1_{(\Si,\cG_{lr},A)})\sseteq\dots$  is the embedding of (germs of) schemes:
\beq
Supp\Big(T^1_{(\Si,\cG_{lr},A)}\Big)\supseteq \cupl^{m-1}_{j=0} \overline{Sing\Big(V\big(I_{j+1}(A)\big)\Big)\smin V\big(I_{j}(A)\big)}.
\eeq
(The closure is taken here in Zariski topology.)

Recall that the singular locus of $V(I_{j+1}(A))$ always contains the locus $V(I_j(A))$. The upper bound in part 3 says that $Supp\Big(T^1_{(\Si,\cG_{lr},A)}\Big)$
 contains the ``unexpected" singular loci of all $V(I_{j+1}(A))$.

The lower bound in part 3 implies in particular that the support of $T^1_{(\Si,\cG_{lr},A)}$ lies inside the
locus $V(I_m(A))$. If $A$ is non-degenerate at some point of $Spec(R)$ then
 $T^1_{(\Si,G,A)}$ is not supported at that point.

\item  This implies: the subscheme  $Supp(T^1_{(\Si,\cG_{lr},A)})\sset Spec(R)$ contains (as subschemes) all
 the components of $V(I_{j+1}(A))\sseteq Spec(R)$ that are not of expected co-dimension. In the classical case, $R=\C\{\ux\}$, the set theoretic version of
  this result is well known, see e.g., \cite{Bruce-Tari04}, \cite{Gusein-Zade-Ebeling}. (And the classical proofs use the
   Thom stratification/transversality theorems, over $\R,\C$.)
\item
\bei
\item Set-theoretically this is equality of the reduced subschemes:
\beq\label{Eq.Results.Support.of.T1.set.theoretic}
Supp\Big(T^1_{(\Si,\cG_{lr},A)}\Big)_{red}=\Big(\cupl^{m-1}_{j=0} \overline{Sing\Big(V\big(I_{j+1}(A)\big)\Big)\smin V\big(I_{j}(A)\big)}\Big)_{red}
 \sset Spec(R).
\eeq
Thus $T^1_{(\Si,\cG_{lr},A)}$ is (set-theoretically) supported exactly on the ``unexpected" singular loci of the determinantal strata.
 This fact is of vital importance and brings numerous corollaries for the determinacy and deformations.
Some versions of this are well known in the classical case, $R=\C\{\ux\}$,  \cite{Bruce-Tari04}, \cite{Gusein-Zade-Ebeling}.
\item The natural wish is to eliminate the radicals in part 5 and to get some bounds close to this equality, e.g., by using the integral closure of
 $Ann(T^1_{(\Si,\cG_{lr},A)})$. This does not seem possible because the cases with
\[
Ann(T^1_{(\Si,\cG_{lr},A)})\sseteq \capl^{m-1}_{j=0}\Big(Sat_{I_{j}(A)}Sing^{(\cm)}_{(m-j)(n-j)}\big(I_{j+1}(A)\big)\Big)^N
\]
  are realized for any $N$, see example \ref{Ex.annihilator.TGlr.cannot.remove.radicals}.
  \eei
\eee

\subsection{The action $\cG_{congr}\circlearrowright  Mat_{m\times m}(R)$, $Mat^{sym}_{m\times m}(R)$,  $Mat^{skew-sym}_{m\times m}(R)$}\label{Sec.Results.Congruence}
We bound the support of $T^1$ in the same way as in theorem \ref{Thm.Results.T1.Grl.coord.changes}.
\bthe\label{Thm.Results.T1.for.Congruence}
Let $R$ be a local Noetherian ring over a field $\k$ of zero characteristic.
\bee[1.]
\item  Let $A\in \Si=Mat_{m\times m}(R)$, $m>1$,  and assume $dim(R)<\lfloor\frac{m}{2}\rfloor$.
 Then $Ann\Big( T^1_{(\Si,\cG_{congr},A)}\Big)\sseteq nil(R)$, the nilradical of $R$.

\item Let $A\in \Si^{sym}:=Mat^{sym}_{m\times m}(R)$.
\bee[i.]
\item
 Suppose a prime ideal $\cp\ssetneq \cm$ satisfies   $ I_{m-1}(A)\not\subseteq \cp\supseteq (det(A))$. Then the localizations at $\cp$ satisfy:
\[
(T^1_{(\Si^{sym},\cG_{congr},A)})_{\cp}=\quot{R_\cp}{\Big((det(A))+Der_\k(R)(det(A))\Big)_\cp}.\]
\item
$Ann.Coker(A)+Ann\quot{\Si^{sym}}{Der_\k(R,\cm)(A)}\sseteq  Ann(T^1_{(\Si^{sym},\cG_{congr},A)})\sseteq
\capl^{m-1}_{j=0}Sat_{I_j(A)}\Big(Sing^{(\cm)}_{\bin{m-j+1}{2}}(I_{j+1}(A))\Big)$.

\item Suppose   for some $j$ the prime decomposition is  $\sqrt{I_{j+1}(A)}=(\cap_\al \cp_\al)\cap(\cap_\be \cq_\be)$,
 where $grade(\cp_\al)<\bin{m-j+1}{2}$ and $grade(\cq_\be)=\bin{m-j+1}{2}$.
 Take the corresponding localizations, $R\stackrel{i_\al}{\to}R_{\cp_\al}$, and suppose
 $\{R_{\cp_\al}\}$ are regular.   Then
$Ann(T^1_{(\cG_{lr},\Si,A)})\sseteq Sat_{I_j(A)}\Big(\cap_\al i^{-1}_\al(I_{j+1}(A)_{\cp_\al})\Big)$.

\item  $\sqrt{Ann\Big(T^1_{(\Si,\cG_{congr},A)}\Big)}=
\capl^{m-1}_{j=0}\Big(\sqrt{Sing^{(\cm)}_{\bin{m-j+1}{2}}(I_{j+1}(A))}:I_j(A)\Big)$.
\eee
\item Let $A\in \Si:=Mat^{skew-sym}_{m\times m}(R)$, $m\ge2$.  Below  $Pf_{m-1}(A)$ is the generalized Pfaffian ideal,   defined in \S\ref{Sec.Background.Determinantal.and.Pfaffians}
\bee[i.]
\item
 Suppose $m$ is even and a prime ideal $\cp\ssetneq\cm$ satisfies $I_{m-2}(A)\not\sseteq \cp\supseteq (det(A))$.
 Then the localizations satisfy
$Ann(T^1_{(\Si,\cG_{congr},A)})_{\cp}=Sing_{1}(I_{m-1}(A))_{\cp}$.
\item
If $m$ is even then
$Ann.Coker(A)+Ann\quots{\Si}{Der_\k(R,\cm)(A)} \sseteq Ann(T^1_{(\Si,\cG_{congr},A)})\sseteq \capl^{m-1}_{\substack{j=0\\j\in2\Z}}Sat_{I_j(A)}\Big( Sing_{\bin{m-j}{2}}\big(I_{j+1}(A)\big)\Big)$.
\item
Suppose $m$ is odd and a prime ideal  $\cp\ssetneq\cm$ satisfies $I_{m-1}(A)\not\sseteq  \cp\supseteq I_{m-3}(A)$.
 Then the localizations satisfy
$Ann(T^1_{(\Si,\cG_{congr},A)})_{\cp}=Sing_{3}(Pf_{m-1}(A))_{\cp}$.

\item If $m$ is odd then
$Pf_{m-1}(A)+Ann\quots{\Si}{Der_\k(R,\cm)(A)} \sseteq Ann(T^1_{(\Si,\cG_{congr},A)})\sseteq \capl^{m-1}_{\substack{j=0\\j\in2\Z}}Sat_{I_j(A)}\Big( Sing_{\bin{m-j}{2}}\big(I_{j+1}(A)\big)\Big)$.

\item   Suppose for some $j$ the prime decomposition   is: $\sqrt{I_{j+1}(A)}=(\cap_\al \cp_\al)\cap(\cap_\be \cq_\be)$,
 where $grade(\cp_\al)<\bin{m-j}{2}$ and $grade(\cq_\be)=\bin{m-j}{2}$.
 Take the corresponding localizations, $R\stackrel{i_\al}{\to}R_{\cp_\al}$ and assume
$\{R_{\cp_\al}\}$ are regular.
  Then
$Ann(T^1_{(\cG_{lr},\Si,A)})\sseteq Sat_{I_j(A)}\Big(\cap_\al i^{-1}_\al(I_{j+1}(A)_{\cp_\al})\Big)$.

\item
For any $m$ (even or odd) holds:
$\sqrt{Ann\Big(T^1_{(\Si,\cG_{congr},A)}\Big)}=
\capl^{m-1}_{\substack{j=0\\j\in2\Z}}\Big(\sqrt{ Sing^{(\cm)}_{\bin{m-j}{2}}\big(I_{j+1}(A)\big)}:I_j(A)\Big)$.
\eee
\eee
\ethe
As in the case of $\cG_{lr}$ the statements have direct geometric interpretations.
\bee[\bf {Part} 1.]
\item   implies: if $dim(R)< \lfloor \frac{m}{2}\rfloor$ then ideal $Ann(T^1)$ contains no power of the maximal ideal.
 This implies (e.g., \cite{Belitski-Kerner}, \cite{Belitski-Kerner2}, \cite{BGK}) that no matrix $A\in Mat_{m\times m}(R)$ is $\cG_{congr}$-finitely determined.
\item Here part i. says that $T^1_{(\Si,\cG_{congr},A)}$ is supported on the essential singular locus of $V(det(A))\sset Spec(R)$,
 with the sublocus $V(I_{m-1}(A))$ erased.
As in the $\cG_{lr}$-case: $A$ is infinitesimally $\cG_{congr}$-stable at the points of $Spec(R)\smin V((det(A)+Der_\k(R)(det(A))))$.

The interpretation of part iii. is as in the same as for $\cG_{lr}$-case.

Part iv. means:  the full (set-theoretic) support of  $T^1_{(\Si,\cG_{congr},A)}$ consists of the ``unexpected" singularities of the
 determinantal strata, cf. equation \eqref{Eq.Results.Support.of.T1.set.theoretic}.
\item
 Part i. can be written also as  $Ann(T^1_{(\Si,\cG_{congr},A)})_{\cp}=Sing^{(\cm)}_{1}(Pf(A))_{\cp}$, using  the properties of
 Pfaffian ideals, \S\ref{Sec.Background.Determinantal.and.Pfaffians}.

The geometric interpretations of the statements are as in Part 2 and for $\cG_{lr}$-case.

The even-odd differences and the  conditions $j\in 2\Z$ are due to peculiarities of $I_j(A)$ for skew-symmetric matrices, e.g.,
 $\sqrt{I_{2j}(A)}=\sqrt{I_{2j-1}(A)}$ and $I_m(A)=0$ for $m$-odd.
\eee

\section{Preparations}\label{Sec.Preparations}
Unless stated otherwise, $R$ is a commutative ring over a field $\k$  of zero characteristic.

\subsection{The module of $\k$-linear derivations}\label{Sec.Background.Der}
 For the (regular) rings $\k[[\ux]]$, $\k\{\ux\}$, $\k\langle\ux\rangle$ one has $Der_\k(R)=R\langle\{\di_j\}\rangle$,
   generated by the first order partial derivatives.

The module of those derivations of $R$ that preserve $\cm$ satisfies: $Der_\k(R,\cm)\supseteq\cm\cdot Der_\k(R)$.
 The equality holds here  for many regular rings.

The module $Der_\k(R)$ localizes nicely, for a prime $\cp\sset R$ holds: $Der_\k(R)_\cp=Der_\k(R_\cp)$, see proposition 16.9 of \cite{Eisenbud}

\subsection{Tangent spaces to the orbits}\label{Sec.Background.Tangent.Spaces}
We recall their presentation e.g., from  \cite[\S3.7]{Belitski-Kerner}.
 The tangent space to the  orbit of a matrix is obtained by applying the tangent space of a group, $T_{(GA,A)} =T_{(G,\one)}A$.
\bee[\bf i.]
\item $G_{lr}=GL(m,R)\times GL(n,R)$ acts by $A\to UAV^{-1}$. Here
\[T_{(G_{lr}A,A)}:=T_{(G_{lr},\one)}A=Mat_{m\times m}(R)\cdot A+A\cdot Mat_{n\times n}(R)\sseteq\Mat.\]
Similarly
 for $G_l$ and $G_r$.
\item  $Aut_\k(R)$. Let $(R,\cm)$ be a local ring, then  $T_{(Aut_\k(R),\one)}=Der_\k(R,\cm)$.
  Here we have only the submodule $Der_\k(R,\cm)\sseteq Der_\k(R)$ because the automorphisms of the local ring
 correspond to the local coordinate changes, i.e.   preserve the origin of $Spec(R)$.
   Therefore
\beq\label{Eq.Tangent.Space.Aut(R)}
 T_{(Aut_\k(R),\one)}A=Der_\k(R,\cm)(A)=\Span_R\{\cD(A)\}_{\cD\in Der_\k(R,\cm)}\sseteq \Mat.
\eeq
\item  $\cG_{lr}:\ A\to U\phi(A)V^{-1}$. Here
$T_{(\cG_{lr},\one)}A=Mat_{m\times m}(R)\cdot A+A\cdot Mat_{n\times n}(R)+Der_\k(R,\cm)(A)\sseteq \Mat$.

\item  $\cG_{congr}:\ A\to U\phi(A)U^t$. Here
$T_{(\cG_{congr},\one)}A=\Span_R\{uA+Au^t\}_{u\in Mat_{m\times m}(R) }+Der_\k(R,\cm)(A)\sseteq Mat_{m\times m}(R)$.
\eee

\subsection{Basic results on localizations}\label{Sec.Background.Localization}
\bel\label{Thm.Background.Localization.1}
Let $\cp\sset R$ be a prime ideal. Then the saturations/localizations satisfy:
\bee[1.]
\item $\cp\not\supseteq J$ iff $J_\cp=R_\cp$.
\item If $\cp\not\supseteq I$ then $Sat_I(J)_\cp=J_\cp$. 
\item Suppose $R$ is Noetherian and fix some ideals $I,J_1,J_2\sset R$. Then
 $Sat_I(J_1)\sseteq Sat_I(J_2)$ iff for any prime $\cp\not\supseteq I$ holds: $(J_1)_\cp\sseteq(J_2)_\cp$.
\eee
\eel
\bpr
\bee[\bf 1.]
\item
See the remark on page 71 of \cite{Bourbaki}.

\item   $\supseteq$ is obvious as $Sat_I(J)\supseteq J$.  For the part $\sseteq$ it is enough to prove:
 $Sat_I(J)\sseteq J_\cp$. Suppose $f\in Sat_I(J)$ then $I^Nf\sseteq J$ for some $N$. As $\cp\not\supseteq I$ we have, by part 1, $I_\cp=R_\cp$.
 Hence $(f)_\cp=(I^Nf)_\cp\sseteq J_\cp$.

\item  $\Rrightarrow$ If   $\cp\not\supseteq I$ then by part (2) holds: $Sat_I(J_i)_\cp=(J_i)_\cp$, thus $(J_1)_\cp\sseteq (J_2)_\cp$.

$\Lleftarrow$ Consider the quotient module
$\quots{Sat_I(J_1)+Sat_I(J_2)}{Sat_I(J_2)}$. The localization of
this quotient at any prime $\cp\not\supseteq I$ vanishes:
 \beq
\Big(\quot{Sat_I(J_1)+Sat_I(J_2)}{Sat_I(J_2)}\Big)_\cp=
\quot{Sat_I(J_1)_\cp+Sat_I(J_2)_\cp}{Sat_I(J_2)_\cp}\stackrel{Part\
2}{=} \quot{(J_1)_\cp+(J_2)_\cp}{(J_2)_\cp}=\{0\}.
 \eeq
 Therefore this quotient is not supported on $Spec(R)\smin V(I)$.
  As $R$ is Noetherian, the ideals are finitely generated and   there exists $N$ satisfying:
 \beq
 Ann(\quot{Sat_I(J_1)+Sat_I(J_2)}{Sat_I(J_2)})\supseteq I^N\quad\Rrightarrow\quad
 I^N\cdot Sat_I(J_1)\sseteq Sat_I(J_2)\quad\Rrightarrow\quad
 Sat_I(J_1)\sseteq Sat_I(J_2).
 \epr \eeq
\eee

{\bf The geometric interpretation.} ($\k=\bk$ is a field.) Take a point $pt\in Spec(R)$ then:
\bee[1.]
\item $pt\not\in V(J)$ iff the germ $(V(J),pt)$ is empty.
\item  Suppose $pt\not\in V(I)$. Then $pt\in\overline{V(J)\smin V(I)}$ iff $pt\in V(J)$.
\item $\overline{V(J_1)\smin V(I)}\supseteq \overline{V(J_2)\smin V(I)}$ iff for any point $pt\not\in V(I)$ holds:
 $(V(J_1),pt)\supseteq (V(J_2),pt)$.
\eee

\bel\label{Thm.Background.Localization.2}
Let $R$ be a local Noetherian ring and  $J_1,J_2\ssetneq R$ some  proper ideals. The following conditions are equivalent:
\bee[1.]
\item $\sqrt{J_1}=\sqrt{J_2}$.
\item For any non-maximal prime ideal, $\cp\ssetneq\cm$, holds: $(J_1)_\cp= R_\cp$ iff $(J_2)_\cp= R_\cp$.
\item  For any non-maximal prime ideal, $\cp\ssetneq\cm$, holds: $(J_1)_\cp\neq R_\cp$ iff $(J_2)_\cp\neq R_\cp$.
\eee
\eel
\bpr
Obviously $2\Leftrightarrow3$, thus we prove $1\Leftrightarrow2$.

$1\Rightarrow2$ If $(J_1)_\cp=R_\cp$ then there exists $f\in J_1$ whose image in $(J_1)_\cp$ is invertible. Thus $f\not\in\cp$. But $f^N\in J_2$ for
 some $N<\infty$. And the image of $f^N$ in $(J_2)_\cp$ is still invertible, hence $(J_2)_\cp=R_\cp$.

$2\Rightarrow1$ Take the prime decomposition, $\sqrt{J_1}=\cap \cp_i$.
(As $\sqrt{J_1}$ is a radical ideal, its primary decomposition consists of  prime ideals.)
Suppose for some $i$ happens $\cp_i\not\supseteq\sqrt{J_2}$ then $\cp_i\not\supseteq J_2$, thus $(J_2)_{\cp_i}=R_{\cp_i}\supsetneq (J_1)_{\cp_i}$. Thus,
 if $\cp$ is a minimal prime for $J_1$ then $\cp_i\supseteq J_2$. Suppose $\cp_i$ is not a minimal prime for $J_2$, then exists a smaller
  prime ideal $\cq\ssetneq\cp_i$, which is a minimal prime for $J_2$. But then, by the same argument as above, $\cq\supseteq J_1$, thus $\cp_i$ could
  not be a minimal prime.

Therefore: $\cp$ is a minimal prime ideal for $\sqrt{J_1}$ iff it is the one for $\sqrt{J_2}$.
In other words, $\sqrt{J_1},\sqrt{J_2}$ have the same minimal primes. Thus, as both are radical, their primary decompositions coincide.
 Hence $\sqrt{J_1}=\sqrt{J_2}$.
\epr

{\em \bf{The geometric interpretation.}} ($\k=\bk$ is a field.) Denote by $0\in Spec(R)$ the base point of the germ. The following are equivalent:
\bee[1.]
\item  Two proper ideals define (set-theoretically) the same locus, $V(J_1)_{red}=V(J_2)_{red}$.
\item For any closed point, $0\neq pt\in Spec(R)$, there holds: $pt\in V(J_1)$ iff $pt\in V(J_2)$.
\eee

\subsection{Saturation vs radicals}
\bel\label{Thm.Background.Saturation}
Given two ideals $I,J\sset R$, with $I$ finitely generated, there holds: $\sqrt{Sat_I(J)}=\sqrt{J}:I$.
\eel
\bpr
$\sseteq$ If $f\in \sqrt{Sat_I(J)}$ then $f^n\in J:I^N$, for some $n,N\in \N$. Thus $f^n\cdot I^N\sseteq J$, hence $f\cdot J\sseteq\sqrt{I}$.

$\supseteq$ Let $I=\langle g_1,\dots,g_n\rangle$ and suppose $f\cdot I\sseteq\sqrt{J}$. Then $f^N\cdot I^N=f^N(g_1,\dots,g_n)^N\sseteq J$ for
 some $N\gg1$.  Thus $f^N\in Sat_I(J)$.
\epr
The finiteness assumption on $I$ is important,   due to the following standard example.
\beq
R=\k[[x_1,x_2\dots]]\supset \cm=\{x_i\}\supseteq J=\langle x_1,x^2_2,x^3_3\dots\rangle.
\quad \text{Then }\sqrt{J}:\cm=R\quad \text{but}\quad\sqrt{Sat_\cm J}=\sqrt{J}=\cm.
\eeq

\subsection{Determinantal and Pfaffian ideals}\label{Sec.Background.Determinantal.and.Pfaffians}
For $1\le j\le m$ and $A\in Mat_{m\times n}(R)$ denote by $I_j(A)\sset R$
the determinantal ideal generated by all the $j\times j$ minors of $A$.
By definition $I_0(A)=R$ and $I_{>m}(A)=\{0\}$.

Determinantal ideals of skew-symmetric matrices, $A\in Mat^{skew-sym}_{m\times m}(R)$, have special properties,
 see e.g., theorem 3.8 in \cite{Ko.La.Sw}. Recall the Pfaffian ideal, $Pf(A)$, and its generalizations $Pf_i(A)=$the ideal generalized by
  Pfaffians of  the principal $i\times i$ submatrices of $A$. We use the following:
\bei
\item For $m$-even: $I_m(A)=Pf(A)^2$ and $I_{m-1}(A)=Pf(A)\cdot Pf_{m-2}(A)$.
\item For $m$-odd: $I_m(A)=0$ and $I_{m-1}(A)=Pf_{m-1}(A)^2$.
\item
 For any $j$ holds:   $\sqrt{I_{2j}(A)}=\sqrt{I_{2j-1}(A)}$.
\eei

\subsection{Annihilator of cokernel}\label{Sec.Background.Ann.Coker} \cite[\S 20]{Eisenbud}
Consider $A\in\Matm$ as a presentation matrix of its cokernel,
$ R^{\oplus n}\stackrel{A}{\to} R^{\oplus m}\to Coker(A)\to0$.

The support of the module $Coker(A)$ is the annihilator-of-cokernel
ideal:
\beq Ann.Coker(A)=Ann\Big(\quot{R^{\oplus m}}{Im(A)}\Big)=\{f\in R|\
f\cdot R^{\oplus m}\sseteq Im(A)\} \sset R.
\eeq
 This ideal is $G_{lr}=GL(m,R)\times GL(n,R)$-invariant and refines the ideal $I_m(A)$.

  The annihilator-of-cokernel is a rather delicate invariant but it is controlled by
 the ideals $\{I_j(A)\}$, see   \cite[proposition 20.7]{Eisenbud} and \cite[exercise 20.6]{Eisenbud}:
\begin{Properties}\label{Eq.Background.Ann.Coker.in.terms.of.Fitting.ideals}
\bee[1.]\item
$\forall\ j<m$:  $Ann.Coker(A)\cdot I_j(A)\sseteq I_{j+1}(A)$,  and
$Ann.Coker(A)^m\sseteq I_m(A)\sseteq Ann.Coker(A)\sseteq\sqrt{I_m(A)}$.
 \item
 If $m=n$ and $det(A)\in R$ is not a zero divisor, then $Ann.Coker(A)=I_m(A):I_{m-1}(A)$.
\item
 If $m<n$ and $grade(I_{m})=(n-m+1)$ then $Ann.Coker(A)=I_m(A)$
\eee\end{Properties}
In particular, for one-row matrices, $m=1$, or when $I_m(A)$ is a radical ideal, $Ann.Coker(A)=I_m(A)$.

\medskip

We use also the following properties of the ideal $Ann.Coker$:
\begin{Properties}\label{Thm.Background.Ann.Coker.Properties}
\bee[1.]
\item  (Block-diagonal case) $Ann.Coker(A\oplus B)=Ann.Coker(A)\cap Ann.Coker(B)$.
\item If $A$ is a square matrix and $det(A)$ is not a zero divisor 
 then $Ann.Coker(A)=Ann.Coker(A^T)$.
\item If  $R$ is a unique factorization domain (UFD) and $A$ is square then $Ann.Coker(A)$ is a principal ideal.
\eee
\end{Properties}
\bpr
\bee[\text{\em Part} 1]
\item  is immediate.

\item  Let $f\in Ann.Coker(A)$ then $AB=f\cdot\one$ for some $B\in Mat_{m\times m}(R)$. Thus $A^\vee\cdot A\cdot B=f\cdot A^\vee$,
 implying: $det(A)B\cdot A=f\cdot det(A)$.
 As $det(A)$ is not a zero divisor we get  $BA=f\cdot \one$, hence $A^tB^t=f\cdot \one$. Thus $f\in Ann.Coker(A^t)$.

\item
 Suppose $0\neq f,g\in Ann.Coker(A)$, then for some $B_f,B_g\in Mat_{m\times m}(R)$ holds: $AB_f=f\cdot\one$, $AB_g=g\cdot\one$.
 Thus (as $g$ is not a zero divisor), by part two we have: $B_gA=g\cdot\one$. Together we get: $g\cdot B_f=B_gAB_f=f\cdot B_g$. Present $g=\tilde{g}\cdot c$, where $(c)\ni f$,
  while $gcd(\tilde{g},f)=1$, i.e., $(\tilde{g})\cap (f)=(\tilde{g}\cdot f)$. Then, as $R$ is UFD, we get: the entries of $B_g$ are divisible by $\tilde{g}$.
   But then $A\cdot(\frac{1}{\tilde{g}}B)=c\cdot\one$, i.e., $c\in Ann.Coker(A)$, and $(c)\ni f,g$.

   Finally, as $R$ is UFD there exists a
    finite decomposition of $c$ into irreducibles. Thus, after a finite such steps we get a generator of $Ann.Coker(A)$.
    \epr
\eee

\beR
\bei
\item Part 2 does not hold when $det(A)$ is nilpotent. For example:
\[
R=\quots{\k[[x,y,z]]}{(y^2,z^2)},\quad A=\bbm x&y\\0&z\ebm,\quad Ann.Coker(A)=(yz,xz),\quad Ann.Coker(A^t)=(xz).
\]
\item
Part 3 does not hold for domains with no unique factorization. For example:
\[
R=\quots{\k[\![x,y,z,w]\!]}{(xy-zw)},\quad  A=\bbm x&0\\0&z\ebm.\quad Ann.Coker(A)=(x)\cap(z)=(xz,xy)\text{ is non-principal.}
\]
\eei
\eeR

\subsection{The generalization of the annihilator of cokernel}\label{Sec.Background.Generaliz.Ann.Coker}
The ideal $Ann.Coker(A)$ is a `partially reduced' version of the ideal of maximal minors $I_m(A)$.
 Equivalently, the annihilator of a module, $Ann(M)$, is a refinement of the minimal Fitting ideal of that module, $Fitt_0(M)$.
More generally, the counterparts of the ideals $\{I_j(A)\}$ (or the Fitting ideals $\{Fitt_{m-j}(M)\}$) are described in \cite{Buchsbaum-Eisenbud}, see also \cite[exercise 20.9]{Eisenbud}.
 We recall briefly the definition and the main properties.

Fix a morphism of free $R$-modules, $E\stackrel{\phi}{\to}F$, here $rank(F)=m<\infty$. For each $1\le j\le m$ define the associated morphism $E\otimes\wedgel^{j-1}F\stackrel{\phi_j}{\to}\wedgel^j F$ by $a\otimes w\to \phi(a)\wedge w$.
\bed
$Ann.Coker_j(\phi):= Ann.Coker(\phi_{m+1-j})$, for $1\le j\le m$.

In addition we define: $Ann.Coker_{j\le 0}(\phi)=R$ and $Ann.Coker_{j>m}(\phi)=0$.
\eed

\begin{Properties}\label{Thm.Background.Ann.Coker.Generalization}
\bee[1.]
\item $Ann.Coker_j(\phi)=Ann\Big(\wedgel^{m+1-j}Coker(\phi)\Big)$. In particular, this ideal is fully determined by the module $Coker(\phi)=\quots{F}{\phi(E)}$.
\item  The ideals $Ann.Coker_j(\phi)$ refine the determinantal ideals, in the following sense:
\bee[i.]
\item
$Ann.Coker(\phi)=Ann.Coker_m(\phi)\sseteq\cdots\sseteq Ann.Coker_j(\phi)\sseteq
 \cdots\sseteq  Ann.Coker_1(\phi)=I_1(\phi)$.
\item For any $i>j\ge0$ holds: $Ann.Coker_i(\phi)\cdot I_j(\phi)\sseteq I_{j+1}(\phi)$.
\item $Ann.Coker_j(\phi)\supseteq I_j(\phi)\supseteq (Ann.Coker_j(\phi))^{j}$.
\item  $I_j(\phi)\sseteq Ann.Coker_j(\phi)\sseteq I_j(\phi):I_{j-1}(\phi)$.   In particular,
 if the image of $I_{j-1}(\phi)$ in $\quots{R}{I_j(\phi)}$
 contains  a non-zero divisor,    then $Ann.Coker_j(\phi)=I_j(\phi)$.
\eee
\item
\bee[i.]
\item  Suppose the map $\phi$ splits block-diagonally, i.e., $E_1\oplus E_2\stackrel{\phi_1\oplus\phi_2}{\to}F_1\oplus F_2$.
Suppose moreover $\phi_1$ is invertible  (thus in particular $rank(E_1)=rank(F_1)$). Then
$Ann.Coker_j(\phi)=Ann.Coker_{j-rank(F_1)}(\phi_2)$.
\item  If $A=diag(\la_1,\dots,\la_m)\in Mat_{m\times m}(R)$ and $(\la_1)\supseteq(\la_2)\supseteq\cdots\supseteq(\la_m)$ then $Ann.Coker_j(\phi)=(\la_j)$.
\eee
\item  The ideals $Ann.Coker_j(\phi)$ are functorial under localizations, i.e., $Ann.Coker_j(\phi)_\cp=Ann.Coker_j(\phi_\cp)$ for any prime $\cp\sset R$.
\item  Suppose $rank(Im(\phi))<r$, then $Ann.Coker_j(\phi)=\{0\}$ for $j\ge r$.
\eee
\end{Properties}
Some remarks/explanations are needed here.
\bee[\bf 1.]
\item Fix some bases of $E,F$, so that $\phi$ is presented by a matrix $A\in \Mat$.  Then
 $Ann.Coker_j(\phi)$ is invariant under $GL(m,R)\times GL(n,R)$-action on $A$. Similarly, fix  $A\in Mat_{m\times n_1}(R)$ and $B\in Mat_{m\times n_2}(R)$.
 If $Span_R(Columns(A))=Span_R(Columns(B))$ then $Ann.Coker_j(A)=Ann.Coker_j(B)$.     If $n_1=n_2$ and  $Span_R(Rows(A))=Span_R(Rows(B))$
  then $Ann.Coker_j(A)=Ann.Coker_j(B)$.
\item 2.i. This sequence of inclusions and the equalities are immediate.
\\2.ii. and 2.iii see  \cite[exercise 20.9]{Eisenbud} and \cite[exercise 20.10]{Eisenbud}.
\\For 2.iv see corollary 1.4. of \cite{Buchsbaum-Eisenbud}.
\item 3.i. In this case $Coker(\phi)\approx Coker(\phi_2)$, now use part 1.
\\3.ii. Follows by explicit check.
\item   Follows straight from $Ann(M_\cp)=Ann(M)_\cp$.
\item Here $rank(Im(\phi))=max\{j|\ I_j(\phi)\neq0\}$.
If $rank(Im(\phi))<r$ then $rank(Im(\phi_{m+1-r}))< rank(\wedgel^{m+1-r}F)$. But then $Ann.Coker(\phi_{m+1-r})=\{0\}$.
\eee

\subsection{The properties of essential singular locus $Sing_r(J)$}\label{Sec.Background.Sing(J)} (defined in \S\ref{Sec.Results.Notations})
First we give an explicit presentation. Let $J=(\uf)$, $Der_\k(R)=\{\cD_\al\}$ be any (not necessarily minimal) choices of  generators.
 Then equation \eqref{Eq.Intro.def.of.Sing(J)} gives:
\beq\label{Eq.Definition.of.Sing(J)}
Sing_r(J) =Ann.Coker_r\bbm \uf&\zero&\zero&\{\cD_\al f_1\}\\\zero&\uf&\zero&\{\cD_\al f_2\}\\\dots&\dots&\dots\\
\zero&\zero&\uf&\{\cD_\al f_N\}\ebm
\eeq
(Here the last column represents the block of columns.)

For the    ideal $Sing^{(\cm)}_r(J)$ one takes $\cD_\al \in Der_\k(R,\cm)$.
\bex
For a principal ideal, $J=(f)$, we get the traditional Tjurina ideal of a function,
\[Sing_1(f)=(f)+Der_\k(R)(f)\sseteq R.\]
More generally, for $J=(f_1,\dots,f_N)$ and $R$ regular,  equation \eqref{Eq.Definition.of.Sing(J)} gives the traditional presentation
 of the singular locus of $V(J)\sset Spec(R)$, but with the annihilator scheme structure instead of the Fitting ideal.
\eex

\subsubsection{Basic properties of $Sing_r(J)$}

Though the definition involved various choices of generators, $Sing_r(J)$ depends on the ideal  $ J$ only.
 Moreover, $Sing_r( J)$  localizes nicely and has other good properties.
\bel\label{Thm.Background.Sing(J).Properties.Basic} Let $R$ be a commutative unital ring.
 Fix an ideal $J \sseteq R$  and some $r\in \N$.
\bee[1.]
\item  The ideal $Sing_r(J)$ does not depend on the choice of the generators of $J$, $Der_\k(R)$.
\item If $J_1 \sset J_2$ then $Sing_r(J_1)\sseteq Sing_r(J_2)$.
\item  $Sing_r(J)\supseteq J+Ann_r\quot{R^{\oplus N}}{Der_\k(R)(\uf)}$, and the inclusion can be proper.
\item  For any prime ideal $\cp\sseteq\cm$ the localization satisfy: $Sing_r(J)_\cp=Sing_r(J_\cp)$. If $\cp\ssetneq\cm$ then
$Sing^{(\cm)}_r(J)_\cp=Sing_r(J_\cp)$.
\item For any  $A\in \Mat$ and any $1\le r,j\le m$ holds: $I_j(A)\sseteq Sing_r(I_j(A))\sseteq I_j(A)+ Der_\k(R)(I_j(A))$.
\eee
\eel
\bpr
\bee[\bf 1.]
\item
Suppose $\uf$, $\tilde\uf$ are two (finite) tuples of generators of $J$. We can assume (extending by zeros) that they are of the same length.
The two tuples are related by $\uf=U\tilde\uf$, $\tilde\uf= V\uf$, for some square  $R$-matrices $U,V$.

Then the matrices in equation  \eqref{Eq.Definition.of.Sing(J)}, for $\uf$, $\tilde\uf$, are related by the left-right multiplication by some $R$-matrices.
  Hence we get: $Sing^{(\uf)}_r( J)\sseteq Sing^{(\tilde\uf)}( J)_r\sseteq Sing^{(\uf)}_r( J)$, and thus
   $Sing^{(\uf)}_r( J)=Sing^{(\tilde\uf)}_r( J)$.

The independence of the choice of generators of $Der_\k(R)$ is even simpler, in this case one should apply only the right multiplications of the matrices.

\item Immediate, just notice that for $\{f_i\}\in J_1$ holds: $\cD_\al(\sum a_i f_i)\equiv \sum a_i \cD_\al(f_i)\ mod(J_1)$.

\item The inclusion is obvious from the presentation in equation \eqref{Eq.Definition.of.Sing(J)},
and the following example shows the possible inequality. Let $\k$ a field of zero characteristic and take $J=(x^7+y^8,x^8+y^9)\sset \k[\![x,y]\!]$.
The height of this ideal is two. We claim that $Sing_2(J)\supseteq(x^8,y^9)$.
Indeed:
\beq
Sing_2(J)=Ann.Coker\bbm x^7+y^8& x^8+y^9&0&0&7x^6&8y^7
\\
0&0&x^7+y^8&x^8+y^9&8x^7&9y^8
\ebm.
\eeq
Denote the columns of this matrix by $\{c_i\}$, then $7c_1+8c_4-x\cdot c_5=(7y^8,8y^9)^t$.
Together with $c_6$ this gives:
\beq
Sing_2(J)=Ann.Coker\bbm
x^7&y^8&0&0 &0&7x^6&8y^7
\\
0&0&x^8&y^9&x^7+y^8&8x^7&9y^8
\ebm.
\eeq
Thus $Sing_2(J)\ni (x^8,y^9)$.  But
\beq
J+Ann_2\quot{R^{\oplus 2}}{Der_\k(R)(\uf)}=(x^7+y^8,x^8+y^9,63x^6y^8-64x^7y^7)=(x^7+y^8,xy^8+y^9,\cm^6\cdot y^{8}).
\eeq
 (The last transition uses the Gr\"{o}bner basis.) From here one sees that e.g.
\beq
 Sing_2(J)\ni x^8\not\in J+Ann_r\quots{R^{\oplus k}}{Der_\k(R)(\uf)}.
\eeq

\item
The equality $Sing_r(J)_\cp=Sing_r( J_\cp)$ holds because the annihilator is functorial on localizations, $Ann_r(M_\cp)=Ann_r(M)_\cp$,
 and the module of derivations  as well (see \S\ref{Sec.Background.Der}).

If $\cp\ssetneq\cm$ then $Der_\k(R,\cm)_\cp\supseteq (\cm\cdot Der_k(R))_\cp=Der_k(R)_\cp=Der_k(R_\cp)$.
Therefore
 $Sing^{(\cm)}_r( J)_\cp=Sing_r( J_\cp)$.

\item
 The inclusion $I_j(A)\sseteq Sing_r(I_j(A))$ holds by part 3. For the inclusion
$Sing_r(I_j(A))\sseteq I_j(A)+ Der_\k(R)(I_j(A))$ we have:
\begin{multline}
Ann_r\quot{R^{\oplus N}}{I_j(A)\cdot R^{\oplus
N}+Der_\k(R)\bbm f_1\\\dots\\ f_N\ebm}\sseteq
\\\sseteq
Ann_r\quot{R^{\oplus N}}{I_j(A)\cdot R^{\oplus
N}+  Der_\k(R)(I_j(A))\cdot R^{\oplus N}}= I_j(A)+  Der_\k(R)(I_j(A)).
\epr
\end{multline}
\eee

\beR
Note that $Der_\k(R)(I_{j+1}(A))\sseteq I_j(A)\cdot Der_\k(R)(I_1(A))$.
 (Expand the $(j+1)\times(j+1)$ minors in terms of $j\times j$ minors.) Therefore
 the upper bound of part 5 of this lemma implies: the singular locus of $V(I_{j+1}(A))$ contains $V(I_j(A))$.
The inclusion  $V\big(Sing(I_{j+1}(A))\big)\supseteq V(I_j(A))$ is often proper.
  However, if $Spec(R)$ is smooth and $A$ is generic then the two sets coincide.
\eeR

\bel\label{Thm.Background.Sing(J).Properties.Advanced}
Suppose $R$ is a local, Noetherian ring.
\bee[1.]
\item  When working with radicals one can replace the annihilator of cokernel by determinantal ideal, $\sqrt{Sing_r(J)}=\sqrt{J+I_r(Der_\k(R)(\uf))}$.
\item Suppose $\cp\sset R$ is a minimal associated prime of $J$ and $grade(\cp)<r$ and $R_\cp$ is a regular ring. Then $Sing_r(J)_\cp=J_\cp$.
\item Suppose the prime decomposition is $\sqrt{J}=(\cap_\al\cp_\al)\cap(\cap_\be \cq_\be)$, where $\{grade(\cp_\al)<r\}$ and $\{grade(\cq_\be)=r\}$.
 Take the corresponding localizations $\{R\stackrel{i_{\cp_\al}}{\to}R_{\cp_\al}\}$ and suppose the rings $\{R_{\cp_\al}\}$ are regular.
   Then
   \[Sing_r(J)\sseteq \capl_{\al}i_\al^{-1}(J_{\cp_\al}) .\]
\eee
\eel
\bpr
\bee[\bf 1.]
\item  By lemma  \ref{Thm.Background.Localization.2} it is enough to verify that for any prime  ideal $\cp\ssetneq\cm$ holds:
\beq
\text{$Sing_r(J)_\cp=R_\cp$ \quad \quad iff \quad \quad $\Big(J+I_r\big(Der_\k(R)(\uf)\big)\Big)_\cp=R_\cp$.}
\eeq
 If $\cp\not\supseteq J$ then both sides are $R_\cp$, as both sides  contain $J$, and $J_\cp=R_\cp$ by lemma \ref{Thm.Background.Localization.1}.
 If $\cp\supseteq J$ then
\beq
\text{$Sing_r(J)_\cp=R_\cp$ \quad  iff \quad
 $\Big(Ann_r\quot{R^{\oplus N}}{Der_\k(R)(\uf)}\Big)_\cp=R_\cp$ \quad  iff \quad
$I_r\Big(Der_\k(R)\big(\uf\big)\Big)_\cp=R_\cp$.}
\eeq

\item
 By part 4  of lemma \ref{Thm.Background.Sing(J).Properties.Basic} we can localize at $\cp$.
 Thus we can assume: $(R,\cm)$ is a  regular local ring and $J\sseteq \cm$.
 Denote by $\ux=(x_1,\dots,x_n)$ a minimal set of generators of the ideal $\cm\sset R$. By the regularity, $n=dim(R)<r$.

Fix some generators $\{f_i\}$ of $J$, we have
\begin{multline}
Sing_r(J)=Ann_r\quot{R^{\oplus N}}{J\cdot R^{\oplus N}+Der_\k(R)\bbm
f_1\\\dots\\f_N\ebm}=\\=Ann_r\quot{\Big(\quots{R}{J}\Big)^{\oplus
N}}{\quots{R}{J}\otimes  Der_\k(R)\bbm f_1\\\dots\\f_N\ebm}=:
Ann.Coker_r^{\quots{R}{J}}\Big[Der_\k(R)[\uf]\Big].
\end{multline}
 Extend the $N$-tuple $(f_1,\dots,f_N)$ to the $N+n$-tuple $(f_1,\dots,f_N,0,\dots,0)$, and compare it to
  the $N+n$-tuple $(f_1,\dots,f_N,x_1,\dots,x_n)$. The latter is a (non-minimal) system of generators of $\cm\sset R$.
 Therefore
\begin{multline}
Ann.Coker_r^{\quots{R}{J}}\Big[Der_\k(R)[\uf]\Big]=
Ann.Coker_r^{\quots{R}{J}}\Big[Der_\k(R)\bbm\uf\\\underline{0}\ebm\Big]\sseteq\\
Ann.Coker_r^{\quots{R}{J}}\Big[Der_\k(R)\bbm\uf\\\ux\ebm\Big]
=
Ann.Coker_r^{\quots{R}{J}}\Big[Der_\k(R)\bbm \ux\ebm\Big].
\end{multline}
Here the two equalities hold by part 1 of lemma \ref{Thm.Background.Sing(J).Properties.Basic}, while the central inclusion holds because $J\sseteq \cm$.

We have obviously
 $Ann.Coker_r^{\quots{R}{J}}\Big[Der_\k(R)\bbm \ux\ebm\Big]\supseteq J$ and it remains to prove the equality.
    As the ring $(R,\cm)$ is Noetherian, the completion is faithful.
  Therefore it is enough to check
  \beq
  \hR\otimes Ann.Coker_r^{\quots{R}{J}}\Big[Der_\k(R)\bbm \ux\ebm\Big]= \hR\cdot J\sset \hR.
\eeq
 By Cohen structure theorem $\hR=\K[[\ux]]$, where $\K\supseteq\k$ is a field.
   Therefore $Der_\k(\hR)=\hR\langle \di_1,\dots,\di_n\rangle+\hR\cdot Der_\k(\K)$. Here $\{\di_i\}$ are the
    classical partial derivatives, while $Der_\k(\K)$ consists of derivations of $\K$, thus $Der_\k(\K)(\ux)=0$.

Therefore we have $Der_\k(\hR)\bbm \ux\ebm=\one_{n\times n}$. Finally, as $n<r$, we get:
\beq
\hR\otimes Ann.Coker_r^{\quots{R}{J}}\Big[Der_\k(R)\bbm \ux\ebm\Big]\sseteq Ann.Coker_r^{\quots{\hR}{J}}\Big[Der_\k(\hR)\bbm \ux\ebm\Big]=\hR\cdot J.
\eeq
Therefore, for the initial ring,
 $Ann.Coker_r^{\quots{R}{J}}\Big[Der_\k(R)\bbm \ux\ebm\Big]=J$.

\item Follows straight from the previous part, just notice $Sing_r(J)\sseteq \cap i^{-1}_\cp Sing_r(J_\cp)$.
\epr
\eee

\beR
\bee[i.]
\item Parts 2,3 read geometrically: the essential singular locus contains all the components of the subscheme $V(J)\sset Spec(R)$
 that are not of expected codimension. (Even if these components are smooth in the classical sense.)
\item The equality in part 3 does not hold, even when $\sqrt{J}=\cp$, with $grade(\cp)<r$.
 For example, let $R=\k[[x,y]]\supset J=(x^p,xy^q)$, of grade $1$.
 Then $Sing_2(J)=(x^p ,xy^q)$. But $\sqrt{J}=(x)$ and $J_{(x)}=(x)$, thus $i^{-1}_{(x)}J_{(x)} =(x)\sset R$.
  \eee
\eeR

\subsubsection{Relation of $Sing_r(J)$ to the classical singular locus of the subscheme $V(J)\sset Spec(R)$}

The singular locus is classically defined using the module of K\"{a}hler
  differentials, $\Om^1_{R/J}$,
 with the Fitting ideal structure,  $Fitt_{dim(\quots{R}{J})} \Om^1_{R/J}\sset R$.

For complete rings in zero characteristic the module of  differentials is often pathological, e.g.
  uncountably generated, see e.g., \S11 of \cite{Kunz}.
  Thus we work in this subsection with universally finite differentials/separated differentials.

For regular rings, the ideal $Sing_r(J)\sset R$
 is the refinement of the classical ideal $Sing(V(J))$, with the annihilator instead of   Fitting scheme  structure:
\bel
Suppose $R$ is a complete regular local Noetherian ring of dimension $n$ and $J\sset R$ is pure of $height=r$.
Then $Sing_r(J)=Ann_{r}\Om^1_{R/J}\sset R$.
\eel
\bpr
For a  complete regular local ring $R$ of dimension $n$,
  and $J=(\uf)\sset R$, the conormal sequence gives, \cite[Proposition 16.3]{Eisenbud}:
\beq
\Om^1_{R/J}=\quots{R}{J}\otimes\frac{R\langle dx_1,\dots,dx_n\rangle }{\{\sum \di_i f_j dx_i\}_{j=1,..,N}}.
\eeq
 As both $Sing_r(J)$ and $Ann_{r}\Om^1_{R/J}$ contain $J$, we compare their images in $\quots{R}{J}$.
 We have the presentation
$(\quots{R}{J})^N\stackrel{\cA}{\to}(\quots{R}{J})^n\stackrel{dx_1,\dots,dx_n}{\to}\Om^1_{R/J}\to0$,
with the presentation matrix
\beq
\cA=
\bbm \di_1 f_1&\dots&\di_1 f_N\\\di_2 f_1&\dots&\di_2 f_N\\\dots&\dots&\dots\\
\di_n f_1&\dots&\di_n f_N\ebm.
\eeq
Here $\cA$  is the transpose of the block of derivatives in equation \eqref{Eq.Definition.of.Sing(J)}.
 Now we notice that $I_r(\cA)$ contains a non-zero divisor modulo $I_{r-1}(\cA)$,
  therefore,
  by part 2.iv of proposition \ref{Thm.Background.Ann.Coker.Generalization},
  \beq
  Fitt_{n-r}(\Om^1_{R/J})=I_r(\cA)=I_r(\cA^t)=Ann_r(\cA^t)=\quots{R}{J}\otimes Sing_r(J).
 \epr\eeq

But in general the two ideals differ essentially, even their radicals differ.
\bex\bee[i.]
\item (The case of non-pure ideal)
Let $R=\k[[x,y,z]]\supset J=(xz,xy)$. Then
\[
Fitt_{dim(R/J)}(\Om^1_{R/J})=Fitt_2\Big(\frac{R\langle dx,dy,dz\rangle}{d(xz),d(xy)}\otimes\quots{R}{J}\Big)=
 I_1 \bbm J&\zero&\zero&z&y\\\zero&J&\zero&0&x\\\zero&\zero&J&x&0\ebm =
 (x,y,z)\sset \quots{\k[[x,y,z]]}{(xz,xy)}.
\]
On the other hand, the expected grade is 2 and :
\[
Sing_2(J)=Ann.Coker_2\bbm J&\zero&z&0&x\\\zero&J&y&x&0\ebm=(x)\sset\k[[x,y,z]].
\]
Thus $\sqrt{Fitt_{2}(\Om^1_{R/J})}\ssetneq \sqrt{Sing_2(J)}$

We observe also: $Sing_1(J)=Ann.Coker_1\bbm J&\zero&z&0&x\\\zero&J&y&x&0\ebm=(x,y,z)\sset\k[[x,y,z]]$.

\item (The case of non-regular rings) Let $R=\quots{\k[[x,y,z]]}{(xy)}\supset J=(z)$. Here $Der_\k(R)=R\langle x\di_x,y\di_y,\di_z\rangle$, thus
 $Sing_1(J)=Ann.Coker_1[0,0,1]=R$. On the other hand:
 \[
\Om^1_{R/J}=\frac{R\langle dx,dy,dz\rangle}{xdy+ydx,dz}\otimes\quots{R}{(z)}\approx \frac{R\langle dx,dy\rangle}{xdy+ydx}\otimes\quots{R}{(z)},\quad
 \text{ thus }  Fitt_1(\Om^1_{R/J})=I_1\bbm z&0&y\\0&z&x\ebm=(x,y,z)\ssetneq Sing_1(J)= R.
 \]
\eee
\eex

\subsection{Invariance of $Ann(T^1)$}\label{Sec.Background.Invariance.Annihilator}
An element $h=(U,V,\phi)\in\cG_{lr}$ acts on $R$ by $f\to\phi(f)$ and $J\to \phi(J)$.

Suppose $h\in \cG_{lr}$ acts on a submodule $\Si\sseteq\Mat$, thus it sends the pair $(\Si,A)$ to the pair $(\Si,hA)$.
\bel\label{Thm.Annihilator.Invariance.Group.Action}
 Let  $G$ be one of the groups $G_l,G_r,G_{lr},Aut_\k(R),\cG_l,\cG_r,\cG_{lr},\cG_{congr}$.
  Suppose $h=(U,V,\phi)\in\cG_{lr}$ acts on $\Si$ and also acts on the $G$-orbits of $A$, i.e., $h(GA)=G(hA)$ .
 Then $\phi(Ann(T^1_{(\Si,G,A)}))=Ann(T^1_{(\Si,G,hA)})$.
\eel
\bpr
 Consider $h$ as a $\k$-linear automorphism of $\Mat$. It induces
 the isomorphism of the tangent
 \\
\parbox{10cm}
{spaces, the first row of the diagram. Its restriction $(\Si,A)\stackrel{}{\to}(\Si,hA)$ induces the second row.
 The restriction $(GA,A)\stackrel{}{\to}(hGA,hA)=(GhA,hA)$ induces the third row.

If $h\in G_{lr}$ then the map $h_*$ is $R$-linear, i.e., $\phi=Id$. If $h\in\cG_{lr}$ then the map is $R$-multiplicative:
}
\quad\quad
$\bM T_{(\Mat,A)}&\stackrel{h_*}{\isom}&T_{(\Mat,hA)}
\\
\cup&&\cup\\
T_{(\Si,A)}&\stackrel{}{\isom}& T_{(\Si,hA)}
\\\cup&&\cup\\
T_{(GA,A)}&\stackrel{}{\isom}&T_{(GhA,hA)}
\eM$

\beq
h_*(f\cdot T_{(\Si,A)})=\phi(f)\cdot h_*(T_{(\Si,A)}),\quad \quad\quad h_*(f\cdot T_{(GA,A)})=\phi(f)\cdot h_*(T_{(GA,A)}).
\eeq
Thus $h$ induces the isomorphism (of $\k$-modules) $T^1_{(\Si,G,A)}\stackrel{h_*}{\to}T^1_{(\Si,G,hA)}$ that satisfies:
$h_*\big(f\cdot T^1_{(\Si,G,A)}\big)=\phi(f)h_*\big(T^1_{(\Si,G,A)}\big)$.
 In particular, if $f\in Ann\big(T^1_{(\Si,G,A)}\big)$ then
 $\phi(f)\in Ann\big(T^1_{(\Si,G,hA)}\big)$, i.e.
 $h^*Ann\big(T^1_{(\Si,G,A)}\big)\sseteq Ann(T^1_{(\Si,G,hA)})$.
As $h$ is invertible, we get the inverse inclusion as well.\epr

\bex
\bee[i.]
\item
 The assumptions of this lemma are obviously satisfied when $h\in G$. This is used to bring $A$ to a particular
 form in the proof of theorems \ref{Thm.Results.T1.Grl.coord.changes}, \ref{Thm.Results.T1.for.Congruence}.
\item In many cases no choice  of $h\in G$ helps,  e.g., $A$ has no nice canonical form under the $G$-action.
Then one takes $h$ in the normalizer of $G$, to ensure $hGA=GhA$. For example, we use the following normal extensions:
 $G_l,G_r\triangleleft G_{lr}$, $\cG_r\triangleleft GL(m,\k)\times\cG_r$, $\cG_l\triangleleft \cG_l\times GL(n,\k)$.
 \item
 Note that $h\in\cG_{lr}\smin\{GL(m,\k)\times\cG_r\}$ does not in general normalize the $\cG_r$-action.
 Similarly for $h\in\cG_{lr}\smin\{\cG_l\times GL(n,\k)\}$.
\eee
\eex

\section{Proofs of the main results}\label{Sec.Proofs}

\subsection{The $\cG_{lr}$-action}\label{Sec.Proofs.cG_lr} \

\bpr (of theorem \ref{Thm.Results.T1.Grl.coord.changes})
\bee[\bf 1.]
\item  Fix some $A=(a_1,\dots,a_n)\in Mat_{1\times n}(R)$. The tangent space $T_{(\cG_{lr}A,A)}$ is written in \S\ref{Sec.Background.Tangent.Spaces}.
 We record the generating matrix of the submodule $T_{(\cG_{lr}A,A)}\sseteq T_{(\Si,A)}= Mat_{1\times n}(R)$:
\beq
\bbm
A&\zero&\dots&\dots&\zero&\{\cD  a_{1}\}\\
\zero&A&\zero&\dots&\zero&\{\cD  a_{2}\}\\
\dots&\dots\\
\zero&\dots&\dots&\zero&A&\{\cD  a_{n}\}
\ebm.
\eeq
(The right column here denotes the block of columns, as $\cD $ runs over the generators of $Der_\k(R,\cm)$.)

Thus $T^1_{(\Si,\cG_{lr},A)}$ is the cokernel of this matrix, while the annihilator of $T^1_{(\Si,\cG_{lr},A)}$, i.e., the $Ann.Coker$ of this matrix,
 is precisely $Sing^{(\cm)}_n\big(I_1(A)\big)$. This proves the equality
 $Ann\big(T^1_{(\Si,\cG_{lr},A)}\big)=Sing^{(\cm)}_{n}(I_1(A))$. The embedding  $Sing^{(\cm)}_{n}(I_1(A))\supseteq I_1(A)+Ann.Coker(Der_\k(R,\cm)(A))$,
  follows by part 3 of lemma \ref{Thm.Background.Sing(J).Properties.Basic} .

\medskip

\item  Fix a prime ideal $\cp$. As $I_{m-1}(A)\not\sseteq\cp$, we get $I_{m-1}(A)_\cp\not\sseteq(\cp)_\cp$ thus $I_{m-1}(A)_\cp= R_\cp$.
 But then at least one of the $(m-1)\times(m-1)$ minors of $A$ becomes invertible in $R_\cp$.
 Therefore the localized matrix is equivalent to a block-diagonal,  $A_\cp\stackrel{(G_{lr})_\cp}{\sim}\one_{(m-1)\times(m-1)}\oplus \tA$.
 As $\cp\supseteq I_m(A)$, we get:  $\tA\in Mat_{1\times (n-m+1)}(\cp_\cp)$.

We assume $A_\cp$ in this form, by \S\ref{Sec.Background.Invariance.Annihilator} such a transition preserves $Ann(T^1_{(\Si,\cG_{lr},A)})_\cp$.
 Then the tangent space to the orbit (see \S\ref{Sec.Background.Tangent.Spaces}) decomposes into the direct sum:
\begin{multline}\label{Eq.inside.proof.decomposition.into.direct.sum}
(T_{(\cG_{lr}A,A)})_\cp= 
 Mat_{m\times m}(R_\cp)\cdot A_\cp+ A_\cp\cdot  Mat_{n\times n}(R_\cp)+Der_\k(R_\cp,\cm_\cp)(A_\cp)=\\
=Mat_{(m-1)\times n}(R_\cp)\oplus Mat_{1\times(m-1)}(R_\cp)\oplus
 \underbrace{\Big(\tA\cdot Mat_{(n-m+1)\times(n-m+1)}(R_\cp)+Der_\k(R_\cp,\cm_\cp)(\tA)\Big)}_{T_{(\cG_{lr}\tA,\tA)}}.
\end{multline}
We use this direct sum decomposition, together with the corresponding direct sum decomposition of $T_{(\Si,A)}$, to get:
\begin{multline}\label{Eq.inside.proof.using.direct.sum}
\Big(Ann(T^1_{(\Si,\cG_{lr},A)}\Big)_\cp=Ann\Big((T^1_{(\Si,\cG_{lr},A)})_\cp\Big)=
Ann\quot{(T_{(\Si,A)})_\cp}{(T_{(\cG_{lr}A,A)})_\cp}\approx\\
\approx Ann\quot{Mat_{1\times
(n-m+1)}(R_\cp)}{T_{(\cG_{lr}\tA,\tA)}}\stackrel{Part\ 1}{=}
Sing^{(\cm_\cp)}_{n-m+1}(I_1(\tA)).
\end{multline}
Note that $I_1(\tA)=I_m(A_\cp)=I_m(A)_\cp$, and the expected height for these ideals is $(n-m+1)$. Altogether:
\beq
Sing^{(\cm_\cp)}_{n-m+1}(I_1(\tA))=Sing^{(\cm_\cp)}_{n-m+1}\big(I_m(A)_\cp\big)=Sing_{n-m+1}\big(I_m(A)\big)_\cp.
\eeq
(For the last step we use $\cp\ssetneq \cm$ and part 4 of lemma \ref{Thm.Background.Sing(J).Properties.Basic}.)

 Thus   $Ann(T^1_{(\Si,\cG_{lr},A)})_{\cp}=Sing_{n-m+1}(I_m(A))_{\cp}$.

\medskip

\item
\bei\item
The embedding
\beq
Ann(T^1_{(\Si,\cG_{lr},A)})\supseteq
Ann.Coker(A)+Ann\quot{\Mat}{Der_\k(R,\cm)(A)}
\eeq
 holds because
$\cG_{lr}\supset G_{r}\rtimes Aut_\k(R)$ gives $T_{(\cG_{lr}A,A)}\supseteq T_{(G_{r}A,A)}+T_{(Aut_\k(R)A,A)}$, and $Ann(T^1_{(\Si,G_{r},A)})=Ann.Coker(A)$.

\item

For any $0\le j<m$ we have to prove: $Ann(T^1_{(\Si,\cG_{lr},A)})\sseteq Sat_{I_{j}(A)}\Big(Sing^{(\cm)}_{(m-j)(n-j)}\big(I_{j+1}(A)\big)\Big)$.

Note that $Ann(T^1_{(\Si,\cG_{lr},A)})\sseteq Sat_{I_j(A)}\Big(Ann(T^1_{(\Si,\cG_{lr},A)})\Big)$, thus it is enough to prove:
\beq
Sat_{I_j(A)}\Big(Ann(T^1_{(\Si,\cG_{lr},A)})\Big)\sseteq Sat_{I_{j}(A)}\Big(Sing^{(\cm)}_{(m-j)(n-j)}\big(I_{j+1}(A)\big)\Big),\quad
 \text{for any } 0\le j\le m-1.
\eeq
Now, by part 3 of lemma \ref{Thm.Background.Localization.1}, it is enough to verify
 the embedding
\beq
Ann(T^1_{(\Si,\cG_{lr},A)})_\cp\sseteq  Sing_{(m-j)(n-j)}\big(I_{j+1}(A)\big)_\cp
\eeq
 for any prime ideal $\cp$ with $\cp\not\supseteq I_{j}(A)$.

Note that $Sing_{(m-j)(n-j)}\big(I_{j+1}(A)\big)\supseteq I_{j+1}(A)$. Thus if $\cp\not\supseteq I_{j+1}(A)$ then
 $Sing_{(m-j)(n-j)}\big(I_{j+1}(A)\big)_\cp=R_\cp$. Therefore it is enough to check only the case when $I_j(A)\not\sseteq\cp\supseteq I_{j+1}(A)$.

 The proof below is similar to that of part two, just for $j<m-1$ we obtain weaker statements.

Take such a prime $\cp$, then  $(I_{j}(A))_\cp=R_\cp$,  by part 1 of lemma \ref{Thm.Background.Localization.1}.
 Thus at least one $j\times j$ minor of $A$ is invertible in $R_\cp$. Therefore the localization of $A$ is block-diagonalizable,
  $A_\cp\stackrel{(G_{lr})_\cp}{\sim} \one_{j\times j}\oplus\tA$, where $\tA\in Mat_{(m-j)\times(n-j)}(R_\cp)$.
By the invariance of annihilator, \S\ref{Sec.Background.Invariance.Annihilator}, we assume $A_\cp$ in this form.

Note that $I_1(\tA)=I_{j+1}(A_\cp)=I_{j+1}(A)_\cp$ and, as $I_{j+1}(A)\sseteq\cp$, we get $I_{j+1}(A)_\cp\sseteq(\cp)_\cp$, i.e., none of the entries of $\tA$ is invertible in $R_\cp$.
As in part two we decompose the tangent space to the orbit into the direct sum.
\begin{multline}\label{Eq.inside.proof.decompos.direct.sum.2}
(T_{(\cG_{lr}A,A)})_\cp=
 Mat_{m\times m}(R_\cp)\cdot A_\cp+A_\cp\cdot  Mat_{n\times n}(R_\cp)+Der_\k(R_\cp,\cm_\cp)(A_\cp)=\\
=Mat_{j\times n}(R_\cp)\oplus Mat_{(m-j)\times j}(R_\cp)\oplus
\underbrace{\Big( Mat_{(m-j)\times(m-j)}(R_\cp)\cdot \tA+\tA\cdot  Mat_{(n-j)\times(n-j)}(R_\cp)+Der_\k(R_\cp,\cm_\cp)(\tA)\Big)}_{T_{(\cG_{lr}\tA,\tA)}}.
\end{multline}
 We simplify the annihilator according to this decomposition:
\begin{multline}
\Big(Ann(T^1_{(\Si,\cG_{lr},A)}\Big)_\cp=Ann\Big((T^1_{(\Si,\cG_{lr},A)})_\cp\Big)=\\=
Ann\quot{(T_{(\Si,A)})_\cp}{(T_{(\cG_{lr}A,A)})_\cp} \approx
Ann\quot{Mat_{(m-j)\times(n-j)}(R_\cp)}{T_{(\cG_{lr}\tA,\tA)}}.
\end{multline}
Unlike part two, for $m-j>1$ we cannot ``pack" the last term in a nice form. Instead we enlarge the annihilator ideal  by observing
that $T_{(G_{lr}\tA,\tA)}\sseteq Mat_{(m-j)\times(n-j)}(I_1(\tA))$. (This is equality for $m-j=1$, but can be a proper embedding for $m-j>1$.)
 Therefore:
\begin{multline}\label{Eq.inside.proof.Annihilator.Localized}
\Big(Ann(T^1_{(\Si,\cG_{lr},A)}\Big)_\cp\sseteq
Ann\quot{Mat_{(m-j)\times(n-j)}(R_\cp)}{Mat_{(m-j)\times(n-j)}(I_1(\tA))+Der_\k(R_\cp,\cm_\cp)(\tA)}=\\=Sing^{(\cm_\cp)}_{(m-j)(n-j)}\big(I_1(\tA)\big).
\end{multline}
Now, as in part two, we observe:
$I_1(\tA)=I_{j+1}(A_\cp)=I_{j+1}(A)_\cp$.
Therefore:
\beq
\Big(Ann(T^1_{(\Si,\cG_{lr},A)}\Big)_\cp\sseteq 
Sing^{(\cm_\cp)}_{(m-j)(n-j)}\big(I_{j+1}(A)_\cp\big)=
\Big(Sing^{(\cm)}_{(m-j)(n-j)}\big(I_{j+1}(A)\big)\Big)_\cp.
\eeq
\eei
As this embedding holds for any localization at $\cp\not\supseteq I_j(A)$, we get (by part 3 of lemma \ref{Thm.Background.Localization.1}):
\beq
Ann(T^1_{(\Si,\cG_{lr},A)})\sseteq Sat_{I_j(A)}\Big(Sing^{(\cm_\cp)}_{(m-j)(n-j)}\big(I_{j+1}(A)\big)\Big).
\eeq

\medskip

\item    Follows right from the bound of part 3 and (part 3 of) lemma \ref{Thm.Background.Sing(J).Properties.Advanced}.

\medskip

\item
 The embedding $\sseteq$ follows by applying lemma  \ref{Thm.Background.Saturation} to the embedding $Ann(T^1_{(\Si,\cG_{lr},A)})\sseteq\dots$ of part three.

  For the embedding $\supseteq$ we use
 lemma \ref{Thm.Background.Localization.2}. Thus it is enough to verify that
 for any non-maximal prime ideal, $\cp\ssetneq\cm$, the localizations of the ideals  satisfy:
\beq
\text{If}\quad
Ann(T^1_{(\Si,\cG_{lr},A)})_\cp\neq R_\cp\quad \text{then} \quad
\capl^{m-1}_{j=0}Sat_{I_j(A)}\Big( Sing_{(m-j)(n-j)}(I_{j+1}(A))\Big)_\cp\neq R_\cp.
\eeq

 Fix a prime ideal $\cp\ssetneq\cm$ and fix the number $r$ satisfying $\cp\supseteq I_{r+1}(A)$, $\cp\not\supseteq I_r(A)$.
  Such $r$ exists, because $I_0(A)=R$, $I_{m+1}(A)=\{0\}$, and is unique, as the chain of ideals $\{I_j(A)\}$ is monotonic.
Moreover, there holds: $0\leq r<m$. If $r=m$ then $I_m(A)\not\sseteq \cp$ implies $I_m(A)_\cp=R_\cp$, hence $Ann(T^1_{(\Si,\cG_{lr},A)})_\cp= R_\cp$.

\

It is enough to prove: $Sat_{I_j(A)}\Big( Sing_{(m-j)(n-j)}(I_{j+1}(A))\Big)_\cp\neq R_\cp$ at least for one value of $j$.
Note that $I_{j+1}(A)_\cp=R_\cp$ for $j+1\le r$, because $I_{r}(A)\not\sseteq\cp$, lemma \ref{Thm.Background.Localization.1}. Therefore
\beq
Sat_{I_j(A)}\Big( Sing_{(m-j)(n-j)}(I_{j+1}(A))\Big)_\cp=R_\cp\quad \text{ for } j+1\le r.
\eeq
Thus we take $j=r$ and prove:
\beq\label{Eq.eq.to.prove}
\text{If}\quad
Ann(T^1_{(\Si,\cG_{lr},A)})_\cp\neq R_\cp\quad \text{then} \quad
Sat_{I_r(A)}\Big( Sing_{(m-r)(n-r)}(I_{r+1}(A))\Big)_\cp\neq R_\cp.
\eeq

As  $\cp\not\supseteq I_r(A)$ we have $I_r(A)_\cp=R_\cp$, so the localization $A_\cp$ of $A$ has at least one invertible minor of size $r\times r$.
 Thus  $A_\cp$ is  $(G_{lr})_\cp$-equivalent to $\one_{r\times r}\oplus\tA$, where $\tA\in Mat_{(m-r)\times(n-r)}(\cp_\cp)$.
 Recall that $Ann(T^1_{(\Si,\cG_{lr},A)})$ is invariant under the $G_{lr}$-equivalence, \S\ref{Sec.Background.Invariance.Annihilator},
 therefore from now on we assume $A_\cp$ in this form.

For this form of $A_\cp$ we have the direct sum decomposition
\beq
(T_{(\cG_{lr}A,A)})_\cp\approx Mat_{r\times n}(R_\cp)\oplus Mat_{(m-r)\times r}(R_\cp)\oplus T_{(\widetilde{\cG_{lr}},\tA)},
\eeq
 as in equations \eqref{Eq.inside.proof.decomposition.into.direct.sum} and \eqref{Eq.inside.proof.decompos.direct.sum.2}. Thus
 $Ann(T^1_{(\Si,\cG_{lr},A)})_\cp\approx Ann(T^1_{(\tilde\Si,\tilde\cG_{lr},\tA)})$, where $\tilde\Si=Mat_{(m-r)\times(n-r)}(R_\cp)$ and $\tilde\cG_{lr}$
  is the corresponding group. Therefore we must prove:
\beq
\text{If}\quad
Ann(T^1_{(\tilde\Si,\tilde\cG_{lr},\tA)})\neq R_\cp\quad \text{then}\quad
 Sat_{I_r(A)}\Big( Sing_{(m-r)(n-r)}(I_{r+1}(A))\Big)_\cp\neq R_\cp.
\eeq
Recall that $T_{(\widetilde{\cG_{lr}},\tA)}=Mat_{(m-r)\times(m-r)}(R_\cp)\cdot \tA+\tA\cdot  Mat_{(n-r)\times(n-r)}(R_\cp)+Der_\k(R_\cp)(\tA)$
  and all the entries of $\tA$ belong to $\cp$.
 Therefore $T_{(\widetilde{\cG_{lr}},\tA)}\ssetneq T_{(\tilde\Si,\tA)}$ iff
$I_1(\tA)\cdot T_{(\tilde\Si,\tA)}+Der_\k(R_\cp)(\tA)\ssetneq T_{(\tilde\Si,\tA)}$.
 Thus $Ann(T^1_{(\tilde\Si,\tilde\cG_{lr},\tA)})\neq R_\cp$ iff $Sing_{(m-r)(n-r)}(I_{1}(\tA))\neq R_\cp$.
Finally, we observe: $Sing_{(m-r)(n-r)}(I_1(\tA))=Sing_{(m-r)(n-r)}(I_{r+1}(A))_\cp$. This proves the implication  \eqref{Eq.eq.to.prove}.

Altogether we have:
\beq
\sqrt{Ann(T^1_{(\Si,\cG_{lr},A)})}=\sqrt{\capl^{m-1}_{j=0}Sat_{I_j(A)}\Big( Sing^{(\cm)}_{(m-j)(n-j)}(I_{j+1}(A))\Big)}.
\eeq
Finally, we apply lemma \ref{Thm.Background.Saturation}.
\epr\eee

\beR\label{Ex.annihilator.TGlr.cannot.remove.radicals}
One would like to lift the radicals in part 5 of this theorem or  to get a better lower bound on $Ann(T^1_{(\Si,\cG_{lr},A)})$ than that of part 3.
This does not seem possible because of the following example. Let $\k$ be a field and  $R$ local Noetherian, $dim(R)=1$.
 (Thus $R$ is the local ring of a germ of curve.) Let $A\in \Matm$, $m>1$, and assume $I_m(A)$ contains a non-zero divisor.
 Then $Sat_{I_j(A)}(J)=R$ for any $j>0$  and any   $ J\sseteq R$ that contains a non-zero divisor. Therefore
 \[
\capl^{m-1}_{j=0}Sat_{I_j(A)}\Big( Sing^{(\cm)}_{(m-j)(n-j)}(I_{j+1}(A))\Big)=Sat_{I_0(A)}\Big( Sing^{(\cm)}_{mn}(I_{1}(A))\Big)= Sing^{(\cm)}_{mn}(I_{1}(A))
\supseteq I_1(A).
 \]
On the other hand, suppose $A$ is diagonal, then: $Ann(T^1_{(\Si,\cG_{lr},A)})=Ann.Coker(A)$.
 Therefore we have trivially
\[
\sqrt{ Ann(T^1_{(\Si,\cG_{lr},A)})}=\sqrt{Ann.Coker(A)}=\cm=\sqrt{I_1(A)}=\sqrt{\capl^{m-1}_{j=0}Sat_{I_j(A)}\Big( Sing^{(\cm)}_{(m-j)(n-j)}(I_{j+1}(A))\Big)}.
\]
But the ideal $I_1(A)$ can be much larger than $Ann.Coker(A)$, and the radicals here cannot be lifted in  any universal way (suitable also for $dim(R)>1$).
\eeR

\subsection{The congruence action, $\cG_{congr}$}\label{Sec.Proofs.Congruence} \

\bpr
 (of theorem \ref{Thm.Results.T1.for.Congruence})
\bee[\bf 1.]
\item
 To verify $Ann(T^1_{(\Si,\cG_{congr},A)})\sseteq nil(R)$ we prove: for any prime  $\cp\sset R$ holds
\beq\label{Eq.vanishing.off.the.origin}
Frac(\quots{R}{\cp})\otimes  Ann(T^1_{(\Si,\cG_{congr},A)})=0\in Frac(\quots{R}{\cp}).
\eeq
This implies: $Ann(T^1_{(\Si,\cG_{congr},A)})$ is inside the intersection of all the prime ideals of $R$, thus inside $nil(R)$.
 (See  e.g., proposition 1.8 of \cite{Atiyah.Macdonald}.)

 Geometrically we check the vanishing of $Ann(T^1_{(\Si,\cG_{congr},A)})$
  at the points of $Spec(R)\smin 0$.

Thus, for any  prime $\cp\sset R$ we study the vector subspace
\begin{multline}
T_{(\cG_{congr}A,A)}\otimes Frac(\quots{R}{\cp})=\Span_{Frac(\quots{R}{\cp})}(UA+AU^t)+Frac(\quots{R}{\cp})\otimes Der_\k(R)(A)\sseteq\\\sseteq
 Mat_{m\times m}\big(Frac(\quots{R}{\cp})\big).
\end{multline}
  As $dim(R)<\lfloor\frac{m}{2}\rfloor$ we have $rank(Der_\k(R)(A))<\lfloor\frac{m}{2}\rfloor$, and therefore
   $dim\Big( Frac(\quots{R}{\cp})\otimes Der_\k(R)(A)\Big)<\lfloor\frac{m}{2}\rfloor$.

  To bound the dimension of $\Span_{Frac(\quots{R}{\cp})}(UA+AU^T)$ we  study  the following vector space of solutions:

\beq
\{U\in Mat_{m\times m}(Frac(\quots{R}{\cp}))|\ UA+AU^T=\zero\}.
\eeq
The later equation is well studied, the dimension of the space of solutions is precisely the codimension of the orbit of $A$ under the congruence.
 The minimal codimension equals
$\lfloor\frac{m}{2}\rfloor$, see e.g., Theorem 3 in \cite{De Teran-Dopico}, and it is achieved for $A\in Mat_{m\times m}(Frac(\quots{R}{\cp}))$ generic. Therefore
 we get:
 \beq
 dim(T_{(G_{congr}A,A)}\otimes Frac(\quots{R}{\cp}))\le m^2- \lfloor\frac{m}{2}\rfloor.
 \eeq
 Therefore $dim(T_{(\cG_{congr}A,A)}\otimes Frac(\quots{R}{\cp}))<m^2$.
 Hence the vanishing in equation \eqref{Eq.vanishing.off.the.origin}.

 Therefore $Ann(T^1_{(\Si,\cG_{congr},A)})$ lies inside the intersection of all the prime ideals $\sseteq \cap\cp$.

\

The proofs of the remaining parts are essentially the same as in theorem \ref{Thm.Results.T1.Grl.coord.changes}, thus we just indicate the main steps.

\item {\bf i.}
Localize at $\cp$, then bring $A_\cp$ to the block-diagonal form $A\stackrel{G_{congr}}{\sim}Diag\oplus \tA$, where
 $Diag\in Mat_{(m-1)\times(m-1)}(R_\cp)$ is invertible, while $\tA\in Mat_{1\times 1}(\cp_\cp)$.
 (See e.g. \cite[theorem 3, page 345]{Birkhoff-MacLane}.)

 Now, as in equation \eqref{Eq.inside.proof.using.direct.sum}, one has
 \beq
 (T^1_{(\Si,\cG_{congr},A)})_\cp\approx\quot{R_\cp}{(\tA)+ Der_\k(R_\cp,\cm_\cp)(\tA)}=\quot{R_\cp}{\Big((det(A))+Der_\k(R,\cm)(det(A)\Big)_\cp}.
 \eeq
Finally, as $\cp\ssetneq \cm$ and $Der_\k(R,\cm)\supseteq \cm\cdot Der_\k(R)$, one has $Der_\k(R,\cm)_\cp=Der_\k(R_\cp)$.


{\bf 2.ii} and {\bf 2.iv}
  The proof is the same as in theorem \ref{Thm.Results.T1.Grl.coord.changes}, just replace $G_{lr}$ by $G_{congr}$.
   We have  $A_\cp\stackrel{G_{congr}}{\sim}Diag\oplus \tA$, where $Diag\in Mat_{j\times j}(R_\cp)$ is invertible.
 Then the analogue of equation \eqref{Eq.inside.proof.Annihilator.Localized} is
\begin{multline}
\Big(Ann(T^1_{(\Si,\cG_{congr},A)}\Big)_\cp\sseteq
Ann\quots{R_\cp}{I_1(\tA)}\otimes\quot{Mat^{sym}_{(m-j)\times(m-j)}(R_\cp)}{Der_\k(R_\cp,\cm_\cp)(\tA)}=
\\=Sing^{(\cm_\cp)}_{\bin{m-j+1}{2}}\big(I_1(\tA)\big)=\Big(Sing^{(\cm)}_{\bin{m-j+1}{2}}\big(I_{j+1}(A)\big)\Big)_\cp.
\end{multline}

{\bf 2.iii}   Follows by  (part 3 of) lemma \ref{Thm.Background.Sing(J).Properties.Advanced}, applied to the bound of part 2.ii.

\item
\bee[\bf i.]
\item ($m$ even) In this case $A_\cp\stackrel{G_{congr}}{\sim}E\oplus\tA$,
 where $E\in Mat^{skew-sym}_{(m-2)\times(m-2)}(R_\cp)$ is invertible, while $\tA\in Mat^{skew-sym}_{2\times 2}(\cp_\cp)$.
 (See e.g.,  \cite[exercise 9, page 347]{Birkhoff-MacLane}.) Thus $I_{m-1}(A)_\cp=I_1(\tA)$.

As before we get:
\beq
(T^1_{(\Si,\cG_{congr},A)})_\cp\approx\quot{Mat^{skew-sym}_{2\times 2}(R_\cp)}{\{u\tA+\tA u^t\}_{u\in Mat_{2\times 2}(R_\cp)}+Der_\k(R,\cm)_\cp(\tA)}.
\eeq
Hence $Ann(T^1_{(\Si,\cG_{congr},A)})_\cp=Sing^{(\cm)}_1(I_1(\tA))=Sing^{(\cm)}_1(I_{m-1}(A))_\cp=Sing_1(I_{m-1}(A))_\cp$. (The last equality because of $\cp\ssetneq\cm.$)

\item ($m$ even)
The bound $\cdots\sseteq Ann(T^1_{(\Si,\cG_{congr},A)})\sseteq \cdots$ is proved as in part 3 of theorem \ref{Thm.Results.T1.Grl.coord.changes}.
 One just notes that $\sqrt{I_{2j}(A)}=\sqrt{I_{2j-1}(A)}$.
\item ($m$ odd)
 In this case  we have $A_\cp\stackrel{G_{congr}}{\sim}E_{(m-3)\times(m-3)}\oplus \tA_{3\times 3}$,
 again by \cite[exercise 9, page 347]{Birkhoff-MacLane}.
 Therefore $I_{m-2}(A)_\cp=I_1(\tA)=Pf_2(\tA)$.

 As before we get:
\beq
(T^1_{(\Si,\cG_{congr},A)})_\cp\approx\quot{Mat^{skew-sym}_{3\times 3}(R_\cp)}{\{u\tA+\tA u^t\}_{u\in Mat_{3\times 3}(R_\cp)}+Der_\k(R,\cm)_\cp(\tA)}.
\eeq
As before, $Der_\k(R,\cm)_\cp=Der_\k(R)$ as $\cp\ssetneq \cm$.

  Recall: $\{u\tA+\tA u^t\}_{u\in Mat_{3\times 3}(R_\cp)}=Mat^{skew-sym}_{3\times 3}(Pf_2(\tA))$, see e.g., lemma 3.5 in \cite{Belitski-Kerner1}.
  Thus
\beq
Ann(T^1_{(\Si,\cG_{congr},A)})_\cp=Sing_3(Pf_2(\tA))=Sing_3(Pf_{m-1}(A)_\cp).
\eeq
\item
For the bound $\cdots\sseteq Ann(T^1_{(\Si,\cG_{congr},A)})$ we note that $T_{(G_{congr},\one)}A\supseteq Mat^{skew-sym}_{m\times m}(Pf_{m-1}(A))$,
see e.g., lemma 3.5. of \cite{Belitski-Kerner}.
 The bound $Ann(T^1_{(\Si,\cG_{congr},A)})\sseteq\cdots$ is proved as before.

\item
As before, this follows by  (part 3 of) lemma \ref{Thm.Background.Sing(J).Properties.Advanced}, applied to the bound of part 3.iv.

\item
As in the proof of theorem \ref{Thm.Results.T1.Grl.coord.changes} we fix a prime ideal
$I_j(A)\not\sseteq\cp\supseteq I_{j+1}(A)$. For $j$ odd one has $\sqrt{I_j(A)}=\sqrt{I_{j+1}(A)}$, (see \S\ref{Sec.Background.Determinantal.and.Pfaffians}),
thus $\cp$ as above exists only for $j$-even. Otherwise the proof is the same.
\epr\eee\eee

\end{document}